\documentclass{amsart}
\usepackage{amssymb}
\usepackage{hyperref}
\usepackage{enumerate}
\theoremstyle{plain}
\newtheorem{theorem}{Theorem}[section]
\newtheorem{lemma}[theorem]{Lemma}
\newtheorem{corollary}[theorem]{Corollary}
\theoremstyle{remark}
\newtheorem{remark}[theorem]{Remark}
\newtheorem{example}[theorem]{Example}

\numberwithin{equation}{section}
\newcommand{\E}{\mathrm{e}}
\newcommand{\I}{\mathrm{i}}
\newcommand{\siul}{\sigma^{\mathrm{u,l}}}
\newcommand{\siu}{\sigma^{\mathrm{u}}}
\newcommand{\sil}{\sigma^{\mathrm{l}}}
\newcommand{\sipmu}{\sigma_\pm^{\mathrm{u}}}
\newcommand{\sipml}{\sigma_\pm^{\mathrm{l}}}
\newcommand{\sipmul}{\sigma_\pm^{\mathrm{u,l}}}
\newcommand{\simpul}{\sigma_\mp^{\mathrm{u,l}}}
\newcommand{\lau}{\lambda^{\mathrm{u}}}
\newcommand{\lal}{\lambda^{\mathrm{l}}}

\renewcommand{\Im}{\mathop{\rm Im}}
\renewcommand{\Re}{\mathop{\rm Re}}
\DeclareMathOperator{\clos}{clos}

\DeclareMathOperator{\loc}{loc}
\DeclareMathOperator{\Res}{Res}
\DeclareMathOperator{\sign}{sign}

\DeclareMathOperator{\inte}{int}
\newcommand{\bal}{\begin{align}}
\newcommand{\eal}{\end{align}}
\newcommand{\nn}{\nonumber}

\newcommand{\de}{\delta}

\newcommand{\pa}{\partial}

\newcommand{\la}{\lambda}
\newcommand{\ov}{\overline}
\def\Xint#1{\mathchoice
   {\XXint\displaystyle\textstyle{#1}}%
   {\XXint\textstyle\scriptstyle{#1}}%
   {\XXint\scriptstyle\scriptscriptstyle{#1}}%
   {\XXint\scriptscriptstyle\scriptscriptstyle{#1}}%
   \!\int}
\def\XXint#1#2#3{{\setbox0=\hbox{$#1{#2#3}{\int}$}
     \vcenter{\hbox{$#2#3$}}\kern-.5\wd0}}
\def\dashint{\Xint-}
\begin{document}
\title[Inverse Scattering for Steplike Finite-Gap Potentials]{Inverse Scattering Theory for One-Dimensional Schr\"odinger Operators with Steplike Finite-Gap Potentials}
\author[A. Boutet de Monvel]{Anne Boutet de Monvel}
\address{Institut de Math\'{e}matiques de Jussieu, case 7012, Universit\'{e} Paris 7\\
2 place Jussieu\\ 75251 Paris \\ France}
\email{\href{mailto:aboutet@math.jussieu.fr}{aboutet@math.jussieu.fr}}

\author[I. Egorova]{Iryna Egorova}
\address{B.Verkin Institute for Low Temperature Physics\\ 47 Lenin Avenue\\61103 Kharkiv\\Ukraine}
\email{\href{mailto:egorova@ilt.kharkov.ua}{egorova@ilt.kharkov.ua}}

\author[G. Teschl]{Gerald Teschl}
\address{Faculty of Mathematics\\
Nordbergstrasse 15\\ 1090 Wien\\ Austria\\ and International Erwin Schr\"odinger
Institute for Mathematical Physics, Boltzmanngasse 9\\ 1090 Wien\\ Austria}
\email{\href{mailto:Gerald.Teschl@univie.ac.at}{Gerald.Teschl@univie.ac.at}}
\urladdr{\href{http://www.mat.univie.ac.at/~gerald/}{http://www.mat.univie.ac.at/\~{}gerald/}}
\thanks{Work supported by the Austrian Science Fund (FWF) under Grant No.\ Y330.
{\it J. d'Analyse Math. {\bf 106:1}, 271--316 (2008)}
}
\begin{abstract}
We develop direct and inverse scattering theory for one-dimensional
Schr\"odinger operators with steplike potentials which are
asymptotically close to different finite-gap  potentials on
different half-axes. We give a complete characterization of the
scattering data, which allow unique solvability of the inverse
scattering problem in the class of perturbations with finite second
moment.
\end{abstract}
\maketitle
\section{Introduction}

In this paper we consider direct and inverse scattering theory for
one-dimensional Schr\"odinger operators with steplike  finite-gap
background, using the approach by Marchenko \cite{M}.

To set the stage, let
\begin{equation}\label{0.1}
H_\pm=-\frac{d^2}{dx^2}+p_\pm(x),\quad x\in \mathbb{R},
\end{equation}
be two (in general different) one-dimensional Schr\"odinger operators
with real  finite-gap potentials $p_\pm(x)$. Furthermore, let
\begin{equation}\label{0.12}
H=- \frac{d^2}{dx^2} +q(x),\quad x\in \mathbb{R},
\end{equation}
be the ``perturbed" operator with real potential $q(x)\in L_{\loc}^1$
such that
\begin{equation}\label{0.2}
\pm \int_0^{\pm \infty}|q(x) - p_\pm(x)| (1+x^2)dx <\infty.
\end{equation}
That is, the potential $q(x)$ has different asymptotic
behavior on different half-axes, and we will call it a
\emph{steplike} potential by analogy with the case of two different
constant backgrounds.

The scattering problem for the operators \eqref{0.12}\textendash\eqref{0.2} is classical and arises in various physical applications, for example, when studying properties of the alloy of two different semi-infinite one dimensional crystals. We refer to the recent work \cite{GNP} for a more detailed discussion of the history of such problems together with further references to the literature. In addition, this scattering problem is of course important for the solution of the Korteweg\textendash de Vries (KdV) equation with initial data in these classes.

For a constant steplike background, that is, $p_\pm(x)=c_\pm$
($c_+\neq c_-$ some constants), this problem was completely solved in
\cite{BF,CK,DS}, including applications to the
initial value problem for the KdV equation (\cite{C,Kap}),
and the asymptotic behavior of the solution for large time
(\cite{Kh}). The case, when the potential vanishes on the right half
axis ($p_+(x)=0$) and is asymptotically periodic finite band on the
left half axis, was considered in \cite{Er,Rb1,Rb2,St}. The initial value problem for the KdV equation was also solved in this case (\cite{Er1}), and long time asymptotics
can be found in \cite{Er1,KhK,KhS} (see also \cite{EBM} for the Jacobi operator case). The case of one-periodic background $p_-(x)=p_+(x)$ was studied in \cite{F1,F2,F3} and
\cite{N1,N2} (for the Jacobi operators the same problem was
considered in \cite{EMT} with the extension to different background
operators in the same isospectral class given in \cite{EMT2}).

However, despite the fact that several special cases are well understood by now, only very little was known about the general situation considered here. In fact, the various mutual locations of the respective background spectra (cf.\ the example on page \pageref{example} below) and background Dirichlet eigenvalues, produce a multitude of different cases. To illustrate this, we mention that only a classification of all possible singularities of the transmission coefficient at the boundary of the spectra would require $16$ case distinctions (this is for example the reason why we formulate its properties in terms of the Wronskian of the Jost solutions in Lemma~\ref{lem2.3}
\textbf{II} below). Our goal here is to find a complete characterization of the scattering data for the operator $H$, that will allow us to solve the inverse scattering problem and to prove the uniqueness of the reconstructed potential in the class \eqref{0.2} with the second moment finite. In particular,
we will do this \emph{without} any restrictions on the mutual location of the respective spectra
and \emph{without} any restrictions on the location of the Dirichlet eigenvalues.
In this respect, note that for example in \cite{F1,F2,F3} it is required that the Dirichlet eigenvalues do not coincide with the edges of the continuous spectrum. In fact, it was demonstrated in \cite{EMT} (for the case of Jacobi operators) that these cases give rise to a different behaviour of the scattering data, which does not occur in the constant background case. Furthermore, the inclusion (and understanding) of this case is important for applications to the solution of the initial value problem of KdV equation, since these cases are unavoidable under the KdV flow. We refer to \cite{EMT1} for a detailed discussion (in case of the Toda lattice) and we will give a brief outline for the KdV equation in Section~\ref{KdV}.

On the other hand, we should also mention that there are two things which we do not
address here: First of all, one could relax our decay assumption and replace the second
moment in \eqref{0.2} by the first moment. Our approach is crafted in such a way that
this causes no principal problems. To keep our presentation more readable we have decidedd
not to include this case at this point. Secondly, one could allow general periodic potentials,
that is, an infinite number of gaps. Again there are no serious impediments to treating
this case.

Finally, let us give a brief outline of the present paper.
We start with some preliminary notations and list some standard facts of the
spectral analysis for the background Hill operators in Section~\ref{secBgOp}.
Then we study the properties of the scattering matrix for steplike
operator, paying particular attention to analytical properties of its entries at the edges
of the continuous spectrum of operator $H$ (Section~\ref{secDSP}).
In Section~\ref{secGLM} we derive the Gel'fand\textendash Levitan\textendash Marchenko (GLM) equations
and obtain complementary estimates on their kernels (see also Appendix~\ref{secTOP}).
In this section we also formulate our main result, that characterize the scattering data
(Theorem~\ref{theor1}).
Then we discuss the unique solvability of the GLM equations, that
allows us to solve the inverse scattering problem. Section~\ref{secISP} is
the most important section of the present paper. Here we discuss
the scheme of the solution of the inverse scattering problem and
prove the uniqueness of the reconstructed potential (Theorem~\ref{theor2}).
Our approach is modeled after the generalized Marchenko approach, developed in \cite{M}.
Our final Section~\ref{KdV} contains some applications to the KdV equation.

\section{The Weyl solutions of the background operators}
\label{secBgOp}

Let $H_\pm$ be two finite-gap one-dimensional Schr\"odinger
operators associated with the potentials
$p_\pm(x)$\footnote{Everywhere in this paper the sub or super index
``$+$" (resp.\ ``$-$") refers to the background on the right (resp.\
left) half-axis.}. Let $s_\pm(z,x)$, $c_\pm(z,x)$ be sin- and
cos-type solutions of the equation
\begin{equation}
\left(-\frac{d^2}{dx^2} +p_\pm(x)\right)y(x)= z\, y(x),
\quad z\in \mathbb{C},\label{1.1}
\end{equation}
associated with the initial conditions
\begin{equation}
s_\pm(z,0)=c^\prime_\pm(z,0) =0,\quad c_\pm(z,0)=s^\prime_\pm(z,0)
=1,\label{1.2}
\end{equation}
where the prime denotes the derivative with respect to $x$.

It is well-known that finite-gap Schr\"odinger operators are
associated with the Riemann surface of a square root of the type
\begin{equation}
\sqrt{-\prod_{j=0}^{2r_\pm} (z-E_j^\pm)}, \qquad E_0^\pm < E_1^\pm <
\cdots < E_{2r_\pm}^\pm,
\end{equation}
where $r_\pm\in \mathbb{N}$. Moreover, $H_\pm$ are uniquely determined by
fixing a Dirichlet divisor
$\sum_{j=1}^{r^\pm}(\mu_j^\pm,\sigma_j^\pm)$, where
$\mu_j^\pm\in[E_{2j-1}^\pm,E_{2j}^\pm]$ and $\sigma_j^\pm\in\{-1, 1\}$.
We refer the interested reader to \cite{BBEIM,GH,GRT,L} for relevant background information.
The reader not familiar with this theory can always think of the special case
of periodic finite gap operators.

The spectra of $H_\pm$ consist of $r_\pm+1$ bands
\begin{equation}
\sigma_\pm:= [E_0^\pm, E_1^\pm]\cup\dots\cup[E_{2j-2}^\pm,
E_{2j-1}^\pm]\cup\dots\cup[E_{2r_\pm}^\pm,\infty). \label{1.61}
\end{equation}
Note that in the special case where $p_\pm$ is periodic, we have
merged all colliding bands. Let
\[
M_{r_\pm} := \left\{\mu_1^\pm, \dots,\mu_{r_\pm}^\pm\right\}
\]
be the set of Dirichlet eigenvalues and set
\begin{equation}\label{1.88}
g_\pm(z)= 
\frac{\prod_{j=1}^{r_\pm}(z -
\mu_j^\pm)}{2\sqrt{-\prod_{j=0}^{2r_\pm} (z -E_j^\pm)}},
\end{equation}
where the branch of the square root is chosen such
that we obtain a Herglotz--Nevanlinna function,
\begin{equation}\label{1.881}
\Im(g_\pm(z)) >0 \quad \mbox{for}\quad \Im(z)>0.
\end{equation}
Let us cut the complex plane along the spectrum $\sigma_\pm$ and denote the upper and
lower sides of the cuts by $\sipmu$ and $\sipml$. The corresponding points on these
cuts will be denoted by $\lau$ and $\lal$, respectively. In particular, this means
\[
f(\lau) := \lim_{\varepsilon\downarrow0} f(\lambda+\I\varepsilon), \qquad
f(\lal) := \lim_{\varepsilon\downarrow0} f(\lambda-\I\varepsilon), \qquad \lambda\in\sigma_\pm.
\]
Condition \eqref{1.881} then implies
\begin{equation}\label{1.8}
\frac{1}{\I} g_\pm(\lau) = \Im(g_\pm(\lau))  >0
\quad \mbox{for}\quad \lambda\in\sigma_\pm.
\end{equation}

Next consider the Weyl solutions $\psi_\pm(z,x)$ and $\breve{\psi}_\pm(z,x)$
of \eqref{1.1} which are  determined up to a multiplication constant, depending on $z$, by the requirement
\begin{equation} \label{1.6}
\begin{split}
&\psi_\pm(z,\,\cdot\,)\in L^2(\mathbb{R}_\pm),\\
&\text{resp. }\breve{\psi}_\pm(z,\,\cdot\,)\in L^2(\mathbb{R}_\mp)
\end{split}
\end{equation}
for $z\in \mathbb{C}\backslash\sigma_\pm$. We will
normalize them according to $\psi_\pm(z,0)=\breve{\psi}_\pm(z,0)=1$ such that
\begin{equation} \label{1.5}
\begin{split}
&\psi_\pm(z,x) = c_\pm(z,x) + m_\pm(z) s_\pm(z,x),\\
&\text{resp. }
\breve{\psi}_\pm(z,x) = c_\pm(z,x) +\breve{m}_\pm(z) s_\pm(z,x),
\end{split}
\end{equation}
where
\begin{equation}
m_\pm(z)= 
\frac{\psi_\pm'(z,0)}{\psi_\pm(z,0)}, \qquad
\breve{m}_\pm(z)= 
\frac{\breve{\psi}_\pm'(z,0)}{\breve{\psi}_\pm(z,0)}, \label{1.4}
\end{equation}
are the Weyl $m$-functions. In the case of periodic operators, $\psi_\pm(z,x)$ and
$\breve{\psi}_\pm(z,x)$ are of course just the Floquet solutions.
They are equal to the branches on the upper/lower sheet of the Baker-Akhiezer function
of $H_\pm$.

It is well-known (see, for example, \cite{L}), that $m_\pm(z)-\breve m_\pm(z)=\mp g_\pm(z)^{-1}$.
Equations \eqref{1.2} and \eqref{1.4} then imply
that the Wronskian of the Weyl solutions is equal to
\begin{equation}\label{1.62}
W(\breve\psi_\pm(z), \psi_\pm(z)) = \mp g_\pm(z)^{-1}.
\end{equation}
where $W(f,g)(x)=f(x)g'(x)-f'(x)g(x)$ denotes the usual Wronski determinant.

The set of band edges is given by
\begin{equation}
\partial\sigma_\pm=\left\{E_0^\pm,E_1^\pm,\dots,E_{2r_\pm}^\pm\right\}.
\end{equation}
For every Dirichlet eigenvalue $\mu^\pm_j$  the Weyl functions might have poles.
If $\mu^\pm_j$ is in the interior of its gap, precisely one Weyl function $m_\pm$ or $\breve m_\pm$
will have a simple pole. Otherwise, if $\mu^\pm_j$ sits at an edge, both will have
a square root singularity. Hence we divide the set of poles accordingly:
\begin{align*}
M_\pm &=\{ \mu^\pm_j\mid\mu^\pm_j \in (E_{2j-1},E_{2j}) \text{ and } m_\pm \text{ has a simple pole}\},\\
\breve M_\pm &=\{ \mu^\pm_j\mid\mu^\pm_j \in (E_{2j-1},E_{2j}) \text{ and } \breve m_\pm \text{ has a simple pole}\},\\
\hat M_\pm &=\{ \mu^\pm_j\mid\mu^\pm_j \in \{E_{2j-1},E_{2j}\} \}.
\end{align*}
Clearly $M_{r_\pm}=M_\pm\cup\breve{M}_\pm\cup\hat M_\pm$.
Then we have
\[
m_\pm(z) =\frac{C_\pm}{z - \mu}\left(1 + o(1)\right),\quad
\breve m_\pm(z)= O(1),
\]
for $z\to \mu\in M_\pm$,
\[
m_\pm(z) =O(1),\quad
\breve{m}_\pm(z)= \frac{\breve C_\pm}{z - \mu}\left(1  + o(1)\right),
\]
for $z\to \mu\in  \breve M_\pm$, and
\[
m_\pm(z) =\frac{C_\pm}{\sqrt{z -E}}\left(1 + o(1)\right),\quad
\breve{m}_\pm(z)=-\frac{ C_\pm}{\sqrt{z -E}}\left(1 + o(1)\right),
\]
for $z\to E\in  \hat M_\pm$. Here $C_\pm$ ,$\breve C_\pm$ denote some nonzero constants.

In particular, we obtain the following properties of the Weyl solutions (see, e.g.,
\cite{GH,GRT,L,M,T}):

\begin{lemma}\label{lem1.1}
The Weyl solutions have the following properties:
\begin{enumerate}[\rm(i)]
\item
The function $\psi_\pm(z,x)$ (resp.\ $\breve\psi_\pm(z,x)$) is holomorphic
as a function of $z$ in the domain $\mathbb{C}\setminus (\sigma_\pm\cup M_\pm)$
(resp.\ $\mathbb{C}\setminus (\sigma_\pm\cup\breve M_\pm)$), takes
real values on the set $\mathbb{R}\setminus \sigma_\pm$, has
simple poles at the points of the set $M_\pm$ (resp., $\breve M_\pm$).
It is continuous up to the boundary $\sipmu\cup \sipml$ except at the points from $\hat M_\pm$ and
\begin{equation}\label{1.10}
\psi_\pm(\lau) = \breve\psi_\pm(\lal) =\overline{\psi_\pm(\lal)},
\quad \lambda\in\sigma_\pm.
\end{equation}
For $E \in \hat M_\pm$ the Weyl solutions satisfy
\[
\psi_\pm(z,x)=O\left(\frac{1}{\sqrt{z-E}}\right), \quad
\breve\psi_\pm(z,x) =O\left(\frac{1}{\sqrt{z-E}}\right), \quad
\mbox{as } z\to E\in \hat M_\pm.
\]
The same is true for $\psi_\pm'(z,x)$ and $\breve \psi_\pm'(z,x)$.
\item
At the edges of the spectrum these functions possess the properties
\[
\psi_\pm(z,x) - \breve\psi_\pm(z,x)=O\left(\sqrt{z -E}\right)\quad \mbox{near}\quad
E\in\pa\sigma_\pm\setminus\hat M_\pm,
\]
and
\[
\psi_\pm(z,x) +
\breve\psi_\pm(z,x)=O(1)\quad\mbox{near}\quad E\in\hat M_\pm\ .
\]
\item
When $z\to\infty$ the following asymptotic behavior holds\/%
\footnote{Here $\Im(\sqrt{z})>0$ as $z\in\mathbb{C}\setminus \mathbb{R}_+$.}:
\[
\psi_\pm(z,x) = \E^{\pm\I\sqrt{z}x}\left(1+O(z^{-1/2})\right)
\quad\mbox{and}\quad
\breve\psi_\pm(z,x) = \E^{\mp\I\sqrt{z}x}\left(1+O(z^{-1/2})\right).
\]
\item
The functions $\psi_\pm(\lambda,x)$ form a complete orthogonal system
on the spectrum with respect to the weight
\begin{equation} \label{1.12}
d\rho_\pm(\lambda) = \frac{1}{2\pi\I} g_\pm(\lambda)d\lambda,
\end{equation}
namely
\begin{equation}\label{1.14}
\oint_{\sigma_\pm}\overline{\psi_\pm(\lambda,y)} \psi_\pm(\lambda,x)d\rho_\pm(\lambda)
= \delta(x-y),
\end{equation}
where $\delta(x)$ is the Dirac delta distribution. Here
we have used the notation
\begin{equation}\label{1.141}
\oint_{\sigma_\pm}f(\lambda)d\rho_\pm(\lambda) := \int_{\sipmu} f(\lambda)d\rho_\pm(\lambda)
- \int_{\sipml} f(\lambda)d\rho_\pm(\lambda).
\end{equation}
\end{enumerate}
\end{lemma}

\section{The direct scattering problem}
\label{secDSP}

Consider the equation
\begin{equation}
\left(-\frac{d^2}{dx^2}+q(x)\right)y(x)= z\, y(x),\quad
z\in \mathbb{C},\label{2.1}
\end{equation}
with a potential $q(x)$, satisfying
condition \eqref{0.2}. This equation has two solutions
$\phi_\pm(z,x)$, the Jost solutions, that are asymptotically close as $x\to
\pm\infty$ to the Weyl solutions of the background equations \eqref{1.1} and
can be represented as (see \cite{F1,F2,F3}):
\begin{equation}
\phi_\pm(z,x)
=\psi_\pm(z,x)\pm\int_{x}^{\pm\infty}\,K_\pm(x,y)\psi_\pm(z,y)\,d y,\label{2.2}
\end{equation}
where $K_\pm(x,y)$ are real-valued, continuously differentiable with respect to
both parameters, and satisfy the estimate
\begin{equation}\label{2.3}
|K_\pm (x,y)|\leq\,\pm C_\pm(x)\,\int_{\frac{x+y}{2}}^{\pm\infty}\,|q(t)-p_\pm(t)|\,d t.
\end{equation}
Here $C_\pm(x)$ are continuous positive functions, monotonically
decreasing (and, therefore, bounded) as $x\to\pm\infty$ (see
Appendix~\ref{secTOP}). For $\lambda\in\sipmu\cup\sipml$ a second pair of solutions of \eqref{2.1}
is given by
\begin{equation}
\overline{\phi_\pm(\lambda,x)
}=\breve\psi_\pm(\lambda,x)\pm\int_{x}^{\pm\infty}\,K_\pm(x,y)\breve\psi_\pm(\lambda,y)
\,d y,\quad\lambda\in\sipmu\cup\sipml.\label{2.4}
\end{equation}
Note $\breve \psi_\pm(\lambda,x) = \ov{\psi_\pm(\lambda,x)}$ for $\lambda\in\sigma_\pm$.

We see that, by formulas \eqref{2.2}, \eqref{2.3}, \eqref{2.4}, and \eqref{1.62},
\begin{equation}\label{2.62}
W\bigl(\phi_\pm(\lambda),\overline{\phi_\pm(\lambda)}\bigr)=\pm g_\pm(\lambda)^{-1}.
\end{equation}
Unlike the Jost solutions, the solutions \eqref{2.4} exist only on the
upper and lower cuts of the spectrum of the corresponding
background, and cannot be continued to the complex plane.

The Jost solutions $\phi_\pm$ are holomorphic in the domains $\mathbb{C}\setminus
\left(\sigma_\pm\cup M_\pm\right)$ and inherit almost all the properties of their
background counterparts, listed in Lemma \ref{lem1.1}, (i)--(ii). In order to remove these
singularities we introduce
\begin{align} \nn
\delta_\pm(z) &:= \prod_{\mu^\pm_j\in M_\pm}(z-\mu^\pm_j),\\ \label{2.6}
\hat \delta_\pm(z) &:= \prod_{\mu^\pm_j\in M_\pm} (z-\mu_j^\pm)
\prod_{\mu^\pm_j \in \hat M_\pm} \sqrt{z - \mu^\pm_j},\\ \nn
\breve \delta_\pm(z) &:= \prod_{\mu^\pm_j\in \breve M_\pm} (z-\mu_j^\pm)
\prod_{\mu^\pm_j \in \hat M_\pm} \sqrt{z - \mu^\pm_j},
\end{align}
where $\prod =1$ if there are no multipliers, and set
\begin{equation}\label{2.12}
\tilde\phi_\pm(z,x)=\delta_\pm(z)
\phi_\pm(z,x),\quad \hat\phi_\pm(z,x)=\hat\delta_\pm(z) \phi_\pm(z,x).
\end{equation}

\begin{lemma}\label{lem2.1}
The Jost solutions $\phi_\pm(z,x)$ have the following properties.
\begin{enumerate}[\rm(i)]
\item
For all $x$, the function $\phi_\pm(z,x)$ considered as
function of $z$, is holomorphic in the domain
$\mathbb{C}\setminus(\sigma_\pm\cup M_\pm)$, takes real values on the set
$\mathbb{R}\setminus\sigma_\pm$, and
has simple poles at the points of the set $M_\pm$. Moreover, $\hat\phi_\pm$ is
continuous up to the boundary $\sipmu\cup \sipml$.
\item
$\phi_\pm(\lambda,x) -\overline{\phi_\pm(\lambda,x)}=O\left(\sqrt{\lambda-
E}\right)$ for
$E\in\pa\sigma_\pm\setminus\hat M_\pm$, and\\
$\phi_\pm(\lambda,x) + \overline{\phi_\pm(\lambda,x)}=O(1)$ for $E\in\hat
M_\pm$.
\end{enumerate}
\end{lemma}

\begin{proof}
Proof of this Lemma follows directly from
\eqref{2.2}, \eqref{2.3}, \eqref{2.4}, Lemma \ref{lem1.1}.
\end{proof}

Now introduce the sets
\begin{equation}\label{2.5}
\sigma^{(2)}:=\sigma_+\cap\sigma_-, \quad
\sigma_\pm^{(1)}=\clos(\sigma_\pm\setminus\sigma^{(2)}), \quad
\sigma:=\sigma_+\cup\sigma_-,
\end{equation}
where $\sigma$ is the (absolutely) continuous spectrum of $H$ and
$\sigma_+^{(1)}\cup\sigma_-^{(1)}$, respectively $\sigma^{(2)}$ are the parts
which are of multiplicity one, respectively two. In addition to the
continuous part, $H$ has a finite number of eigenvalues situated
in the gaps, $\sigma_d\subset\mathbb{R}\setminus\sigma$ (see, e.g., \cite{R-B}).
We will use the notation $\inte(\sigma_\pm)$ for the interior of the spectrum,
that is, $\inte(\sigma_\pm):=\sigma_\pm\setminus\pa\sigma_\pm$.

The following example illustrates the various possible locations of
the spectra together with the Dirichlet eigenvalues.
\vskip 0.2cm

\begin{example}                 \label{example}
Let $H_+$ be the two-band quasi-periodic operator with the spectrum on
the set $\sigma_+=[E_1, E_2]\cup[E_4, +\infty)$ and $H_-$ be the
three band operator with the spectrum $\sigma_-=[E_1, E_2]\cup[E_3,
E_4]\cup [E_5,+\infty)$, where $E_1<E_2<\cdots<E_5$ (cf.\ Figure~\ref{figsi}).
Then $\sigma=[E_1, E_2]\cup[E_3, +\infty)$, $\sigma_+^{(1)} =[E_4,E_5]$,
$\sigma_-^{(1)}=[E_3,E_4]$, and $\sigma^{(2)}=[E_1,E_2]\cup[E_5,+\infty)$.
Denote by $\mu^-_1$ the Dirichlet eigenvalue for the operator $H_-$,
that belongs to the closed gap $[E_2, E_3]$ and let $\mu^+_1$ be the
Dirichlet eigenvalue of $H_+$ from the gap $[E_2,E_4]$.
\begin{figure}[ht]
\begin{picture}(11,1.2)
\put(1,0.2){$\sigma_-$}
\put(1,0.7){$\sigma_+$}

\put(1,0.5){\line(1,0){8}}
\put(3,0.4){\rule{9mm}{2mm}}
\put(2.9,0){$\scriptstyle E_1$}
\put(3.8,0){$\scriptstyle E_2$}
\put(4.4,0.42){$\bullet$}
\put(4.3,0.1){$\scriptstyle \mu^-_1$}
\put(5,0.4){\rule{11mm}{1mm}}
\put(4.9,0){$\scriptstyle E_3$}
\put(6,0){$\scriptstyle E_4$}
\put(5.3,0.42){$\bullet$}
\put(5.2,0.8){$\scriptstyle \mu^+_1$}
\put(6.1,0.5){\rule{39mm}{1mm}}
\put(6.8,0.42){$\bullet$}
\put(6.7,0.1){$\scriptstyle \mu^-_2$}
\put(7.7,0.4){\rule{23mm}{1mm}}
\put(7.6,0){$\scriptstyle E_5$}
\end{picture}
\caption{Typical mutual locations of $\sigma_-$ and $\sigma_+$.}\label{figsi}
\end{figure}
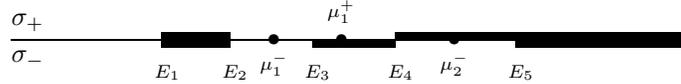
\end{example}

\noindent
We may encounter various mutual locations of these eigenvalues. For example:
\begin{enumerate}[\textbullet]
\item
$\mu_1^+\neq\mu_1^-$ and $\mu_1^+, \mu_1^-\in (E_2, E_3)$ (i.e., $\mu_1^\pm\in M_\pm$),
\item
$\mu_1^+=\mu_1^-\in (E_2, E_3)$,
\item
$\mu_1^+=E_2, \mu_1^-\neq E_2$ (i.e. $\mu_1^+\in \hat M_+$),
\item
$\mu_1^+\in[E_3, E_4]$ (the Dirichlet eigenvalue is situated on the
spectrum of multiplicity one) and $\mu_1^-\neq E_3$,
\item
etc.
\end{enumerate}
All such mutual locations of the Dirichlet eigenvalues imply
different properties of the scattering data and have to be studied
separately.

Let
\begin{equation}\label{2.8}
W(z):=W(\phi_-(z,\,\cdot\,), \phi_+(z,\,\cdot\,))
\end{equation}
be the Wronskian of two Jost solutions. This function is meromorphic in
the domain $\mathbb{C}\setminus\sigma$ with possible poles at the
points $M_+\cup M_-\cup(\hat M_+\cap\hat M_-)$
and with possible square root singularities at the points $\hat
M_+\cup\hat M_-\setminus (\hat M_+\cap\hat M_-)$. Set
\begin{equation}\label{2.10}
\tilde W(z)=W(\tilde\phi_-(z),\tilde\phi_+(z)),
\quad\hat W(z)=W(\hat\phi_-(z),\hat\phi_+(z)).
\end{equation}
Then the function $\hat W(\lambda)$ is holomorphic in
the domain $\mathbb{C}\setminus\mathbb{R}$
and continuous up to the boundary. But unlike the functions $W(z)$
and $\tilde W(z)$, it may not take real values on the set
$\mathbb{R}\setminus\sigma$ and complex conjugated values on the
different sides of the spectrum. That is why it is more convenient
to characterize the spectral properties of the operator $H$ by means
of the function $\tilde W$, which can have singularities at the
points of the set $\hat M_+\cup\hat M_-$. We will study the precise
character of these singularities below.

Since the discrete spectrum of our operator $H$ is finite, we can write it as
\[
\sigma_d=\{\lambda_1,\dots,\lambda_p\}\subset\mathbb{R}\setminus\sigma.
\]
For every eigenvalue
we introduce the corresponding
norming constants
\begin{equation} \label{2.14}
\left(\gamma_k^\pm\right)^{-2}=\int_{\mathbb{R}}\,\tilde\phi_\pm^2(\lambda_k,x)\,dx.
\end{equation}

Note that outside the spectrum $\tilde W(z)=0$ vanishes precisely at the
eigenvalues. However, it might also vanish inside the spectrum at points in $\pa\sigma_-\cup\pa\sigma_+$
and we will call such points virtual levels of the operator $H$
\begin{equation}\label{2.363}
\sigma_v:=\{E\in\sigma\mid\hat W(E)=0\}.
\end{equation}
We will show $\sigma_v \subseteq
\pa\sigma\cup(\pa\sigma_+^{(1)}\cap\pa\sigma_-^{(1)})$ in Lemma~\ref{lem2.3} below.
All other points $E$ of the set $\pa\sigma_+\cup\pa\sigma_-$ correspond to
the generic case $\hat W(E)\neq 0$.

Now we begin our study of the properties of the scattering matrix. Introduce the
scattering relations
\begin{equation}\label{2.16}
T_\mp(\lambda) \phi_\pm(\lambda,x)
=\overline{\phi_\mp(\lambda,x)} + R_\mp(\lambda)\phi_\mp(\lambda,x),
\quad\lambda\in\simpul,
\end{equation}
where the transmission and reflection coefficients are defined as usual,
\begin{equation}\label{2.17}
T_\pm(\lambda):= \frac{W(\overline{\phi_\pm(\lambda)}, \phi_\pm(\lambda))}{W(\phi_\mp(\lambda), \phi_\pm(\lambda))},\qquad
R_\pm(\lambda):= - \frac{W(\phi_\mp(\lambda),\overline{\phi_\pm(\lambda)})}{W(\phi_\mp(\lambda),\,\phi_\pm(\lambda))}, \quad\lambda\in \sipmul.
\end{equation}
Their characteristic properties will be given in the following lemma.

\begin{lemma}\label{lem2.3}
For the entries of the scattering matrix
the following properties are valid:
\begin{enumerate}[\bf I.]
\item
\begin{enumerate}[\bf(a)]
\item
$T_\pm(\lau) =\overline{T_\pm(\lal)}$ and
$R_\pm(\lau) =\overline{R_\pm(\lal)}$ for $\lambda\in\sigma_\pm$.
\item
$\dfrac{T_\pm(\lambda)}{\overline{T_\pm(\lambda)}}= R_\pm(\lambda)$
for $\lambda\in\sigma_\pm^{(1)}$.
\item
$1 - |R_\pm(\lambda)|^2 =
\dfrac{g_\pm(\lambda)}{g_\mp(\lambda)}\,|T_\pm(\lambda)|^2$ for $\lambda\in\sigma^{(2)}$.
\item
$T_\pm(\lambda) = 1 + O(|\lambda|^{-1/2})$ and
$R_\pm(\lambda) = O(|\lambda|^{-1/2})$ for $\lambda\to\infty$.
\item
$\overline{R_\pm(\lambda)}T_\pm(\lambda) +
R_\mp(\lambda)\overline{T_\pm(\lambda)}=0$ for $\lambda\in\sigma^{(2)}$.
\end{enumerate}
\item
The functions $T_\pm(\lambda)$ can be extended analytically to the domain
$\mathbb{C} \setminus (\sigma\cup M_\pm\cup\breve M_\pm)$ and satisfy
\begin{equation}\label{2.18}
\frac{-1}{T_+(z) g_+(z)} = \frac{-1}{T_-(z) g_-(z)}=:W(z),
\end{equation}
where $W(z)$ possesses the following properties:
\begin{enumerate}[\bf(a)]
\item
The function $\tilde W(z)=\delta_+(z)\delta_-(z) W(z)$ is holomorphic in the
domain $\mathbb{C}\setminus\sigma$, with simple zeros at the points $\lambda_k$, where
\begin{equation}\label{2.11}
\biggl(\frac{d\tilde W}{d z}(\lambda_k)\biggr)^2
=(\gamma_k^+\gamma_k^-)^{-2}.
\end{equation}
Besides, it satisfies
\begin{equation}\label{2.9}
\overline{\tilde W(\lau)}=\tilde W(\lal),
\quad \lambda\in\sigma\quad \text{and}\quad \tilde W(\lambda)\in\mathbb{R}
\quad \text{for} \quad \lambda\in\mathbb{R}\setminus \sigma.
\end{equation}
\item
The function $\hat W(z) = \hat\delta_+(z) \hat\delta_-(z) W(z)$ is
continuous on the set $\mathbb{C}\setminus\sigma$ up to the boundary
$\siu\cup\sil$. It can have zeros on the set
$\pa\sigma\cup(\pa\sigma_+^{(1)}\cap\pa\sigma_-^{(1)})$ and does not vanish
at the other points of the set $\sigma$.  If $\hat W(E)=0$ as
$E\in\pa\sigma\cup(\pa\sigma_+^{(1)}\cap\pa\sigma_-^{(1)})$, then $\hat W(z)
=  \sqrt{z -E} (C(E)+o(1))$, $C(E)\ne 0$.
\end{enumerate}
\item
\begin{enumerate}[\bf(a)]
\item
The reflection coefficient $R_\pm(\lambda)$, is a continuous function on the set
$\inte(\sipmul)$.
\item
If $E\in\pa\sigma_+\cap \partial\sigma_-$ and $E\notin \sigma_v$ then the function $R_\pm(\lambda)$
is also continuous at $E$. Moreover
\begin{equation}\label{P.1}
R_\pm(E)=
\begin{cases}
-1 &\text{for } E\notin\hat M_\pm,\\
1 &\text{for } E\in\hat M_\pm.
\end{cases}
\end{equation}
\end{enumerate}
\end{enumerate}
\end{lemma}

\begin{proof}
The  proof  is based on formulas \eqref{2.17}, \eqref{2.62}, \eqref{2.2}, and
Lemma~\ref{lem2.1}.

\textbf{I}.
The symmetry property \textbf{(a)} follows from formulas \eqref{2.17}, \eqref{1.10},
\eqref{2.2}, and \eqref{2.4}. To verify \textbf{(b)} observe, that
$\tilde\phi_\mp(\lambda,x)\in\mathbb{R}$ as
$\lambda\in\inte(\sigma_\pm^{(1)})$. Together with \eqref{2.17}
and \eqref{2.12} this implies {\bf(b)}. Now let $\lambda\in\inte(\sigma^{(2)})$.
Then by \eqref{2.16}
\[
|T_\pm|^2W(\phi_\mp,\overline{\phi_\mp}) = (|R_\pm|^2 -1)
W(\phi_\pm,\overline{\phi_\pm})
\]
and property \textbf{(c)} follows from \eqref{2.62} and \eqref{1.12}. To prove
\textbf{(d)} we use \eqref{2.62}, \eqref{1.8} and \eqref{2.2}. Then
property (iii) of Lemma~\ref{lem1.1} allows us to proceed as in
the proof of \cite[Lemma~3.5.1]{M}to obtain the necessary
asymptotics. The consistency condition \textbf{(e)} and the identity
\eqref{2.18}, considered on $\inte(\sigma^{(2)})$, can be
derived directly from the definition \eqref{2.17}.

\textbf{II}.
\textbf{(a)}. Except for \eqref{2.11} everything follows from the corresponding
properties of $\phi_\pm(z,x)$ and it remains to show \eqref{2.11}.
If $\hat W(\lambda_0)=0$ for some $\lambda_0\in\mathbb{C}\setminus\sigma$, then
\begin{equation}\label{D.1}
\tilde\phi_\pm(\lambda_0,x)=c_\pm\tilde\phi_\mp(\lambda_0,x)
\end{equation}
for some constants $c_\pm$ (depending on $\lambda_0$) and satisfying $c_-\,c_+ =1$.
In particular, each zero of $\tilde W$ (or $\hat W$) outside the continuous spectrum,
is a point of the discrete spectrum of $H$ and vice versa.

Let $\gamma_\pm$ be the norming constants defined in \eqref{2.14} for some point of the discrete spectrum
$\lambda_0$. Proceeding as in \cite{M} one obtains
\begin{equation}\label{D.2}
W\Bigl(\tilde\phi_\pm(\lambda_0,0),\frac{d}{d\lambda}\tilde\phi_\pm(\lambda_0,0)\Bigr)=
\int_0^{\pm\infty}\tilde\phi_\pm^2(\lambda_0,x)\,dx.
\end{equation}
Equalities \eqref{D.1} and \eqref{D.2} imply
\begin{align} \nn
\gamma_\pm^{-2} &=
\mp c_\pm^2\int_0^{\mp\infty} \tilde\phi_\mp^2(\lambda_0,x) dx \pm \int_0^{\pm\infty}
\tilde\phi_\pm^2(\lambda_0,x) dx\\ \nn
& =\mp\ c_\pm^2W\Bigl(\tilde\phi_\mp(\lambda_0,0),\frac{d}{d\lambda}\tilde\phi_\mp(\lambda_0,0)\Bigr) \pm
W\Bigl(\tilde\phi_\pm(\lambda_0,0),\frac{d}{d\lambda}\tilde\phi_\pm(\lambda_0,0)\Bigr)\\ \label{D.3}
&= c_\pm\,\frac{d}{d\lambda}W(\tilde\phi_-(\lambda_0),\tilde\phi_+(\lambda_0))
\end{align}
and, since $c_-c_+=1$, we obtain \eqref{2.11}.

Item \textbf{(b)} will be shown in Lemma~\ref{lem2.2}.

\textbf{III},
\textbf{(a)} follows from the corresponding properties of $\phi_\pm(z,x)$ and from \textbf{II, (b)}.
To show \textbf{III, (b)} we use that by \eqref{2.17} the reflection coefficient has
the representation
\begin{equation}\label{2.n}
R_\pm(\lambda)=-\frac{W(\ov{\phi_\pm(\lambda)},\phi_\mp(\lambda))}
{W\left(\phi_\pm(\lambda),\phi_\mp(\lambda)\right)}=\pm\frac{W(\ov{\phi_\pm(\lambda)},
\phi_\mp(\lambda))}{W(\lambda)}
\end{equation}
and is continuous on both sides of the set $\inte(\sigma_\pm)\setminus (M_\mp\cup\hat M_\mp)$.
Moreover,
\[
|R_\pm(\lambda)|=\left|\frac{W(\ov{\hat\phi_\pm(\lambda)},\hat\phi_\mp(\lambda))}{\hat W(\lambda)}\right|,
\]
where the denominator does not vanish on the set $\sigma_\pm\setminus\sigma_v$. Hence
$R_\pm(\lambda)$ is continuous on this set since both the numerator and denominator are.

Next, let $E\in\pa\sigma_\pm \setminus \sigma_v$ (in particular $\hat
W(E)\ne 0$). Then, if $E\notin \hat M_\pm$, we use \eqref{2.n} in the form
\begin{equation}\label{2.n1}
R_\pm(\lambda)= -1 \mp \frac{\hat{\de}_\pm(\lambda)
W(\phi_\pm(\lambda)-\ov{\phi_\pm(\lambda)},\hat\phi_\mp(\lambda))}{\hat W(\lambda)},
\end{equation}
which shows $R_\pm(\lambda) \to -1$ since $\phi_\pm(\lambda)-\ov{\phi_\pm(\lambda)}\to 0$ by Lemma~\ref{lem2.1} (ii). This settles the first case in \eqref{P.1}. Similarly, if
$E\in\hat M_\pm$, we use \eqref{2.n} in the form
\begin{equation}\label{2.n2}
R_\pm(\lambda)= 1 \pm \frac{\hat{\de}_\pm(\lambda)
W(\phi_\pm(\lambda)+\ov{\phi_\pm(\lambda)},\hat\phi_\mp(\lambda))} {\hat W(\lambda)},
\end{equation}
which shows $R_\pm(\lambda) \to 1$ since
$\hat{\de}_\pm(\lambda)\to 0$ and $\phi_\pm(\lambda)+\ov{\phi_\pm(\lambda)} =
O(1)$ by Lemma~\ref{lem2.1} (ii). This settles the second case in
\eqref{P.1} as well.
\end{proof}

We note that the behavior of $T_\pm(z)$ near the boundary points of the spectra
can be read off from
\begin{equation}\label{2.19}
T_\pm(z) = \frac{-1}{g_\pm(z) W(z)} = -
\frac{\hat\delta_\mp(z)}{\breve\delta_\pm(z)}
\frac{2 \sqrt{\prod_{j=0}^{2r_\pm} (z-E_j^\pm)}}{\hat W(z)}.
\end{equation}

\section{The Gel'fand-Levitan-Marchenko equation}
\label{secGLM}

The aim of this section is to derive the inverse scattering
problem equation (the Gel'fand-Levitan-Marchenko equation) and to discuss some
additional properties of the scattering data, that are consequences of this equation.

\begin{lemma}\label{lem4.2}
The inverse scattering problem (the GLM) equation has the form
\begin{equation}\label{ME}
K_\pm(x,y) + F_\pm(x,y) \pm \int_x^{\pm\infty} K_\pm(x,t) F_\pm(t,y)d t =0,
\quad \pm y>\pm x
\end{equation}
where
\begin{align}\label{4.2}
F_\pm(x,y)
&= \oint_{\sigma_\pm} R_\pm(\lambda)\,\psi_\pm(\lambda,x) \psi_\pm(\lambda,y)
d\rho_\pm(\lambda) + \\ \nn
&\quad
+\,\int_{\sigma_\mp^{(1),\mathrm{u}}} |T_\mp(\lambda)|^2 \psi_\pm(\lambda,x)
\psi_\pm(\lambda,y)d\rho_\mp(\lambda)\\ \nn
&\quad
+ \sum_{k=1}^p (\gamma_k^\pm)^2 \tilde\psi_\pm(\lambda_k,x) \tilde\psi_\pm(\lambda_k,y).
\end{align}
\end{lemma}

\begin{proof}
Consider the function
\begin{align}\label{4.7}
G_\pm(z,x,y) &=
T_\pm(z) \phi_\mp(z,x)\psi_\pm(z,y) g_\pm(z)
- \breve\psi_\pm(z,x)\psi_\pm(z,y) g_\pm(z)\\ \nn
& :=G_\pm^\prime(z, x,y)+G_\pm^{\prime\prime}(z, x,y),\quad \pm y>\pm x,
\end{align}
where $x,y$ are considered as parameters.  As a function of $z$ it is
meromorphic in the domain $\mathbb{C}\setminus\sigma$ with simple poles
at the points $\lambda_k$ of the discrete spectrum. It is continuous up
to the boundary $\siu\cup\sil$, except for the points of the edges of background
spectra,  where
\begin{equation}\label{4.8}
G_\pm(z,x,y)=O\bigl((z -E)^{-1/2}\bigr) \quad \text{as}\quad
E\in\pa\sigma_+\cup\pa\sigma_-.
\end{equation}
Since for $z\to\infty$ we have
\begin{align*}
\phi_\mp(z,x) & = \E^{\mp \I\sqrt z x}\bigl(1+O(z^{-1/2})\bigr),\quad
\breve \psi_\mp(z,x) = \E^{\mp \I\sqrt z x}\bigl(1+O(z^{-1/2})\bigr),\\
\psi_\pm(z,y) &= \E^{\pm \I\sqrt z y}\bigl(1+O(z^{-1/2})\bigr),\quad
T_\pm(z)=1+O(z^{-1/2})
\end{align*}
and $g_\pm(z) = \frac{-1}{2\I\sqrt z} +O(z^{-1})$ then
\begin{equation}\label{4.9}
G_\pm(z,x,y)= \E^{\pm \I\sqrt z \,(y-x)} O(z^{-1}), \quad \pm y>\pm x.
\end{equation}
Consider a closed contour $\Gamma_\varepsilon$ consisting of a large circular arc
together with some parts wrapping around the spectrum $\sigma$ inside this
arc at a small distance from the spectrum. By the Cauchy theorem
\[
\frac{1}{2\pi\I}\oint_{\Gamma_\varepsilon}\ G_\pm(z,x,y)d z =
\sum_{\lambda_k\in\sigma_d}\Res_{\lambda_k}G_\pm(z,x,y).
\]
Estimate \eqref{4.9} allows us to apply Jordan's lemma. Since  by
\eqref{4.8} the limit value of $G_\pm$ as $\varepsilon\to 0$ is
integrable on $\sigma$, and the function $G^{\prime\prime}_\pm$ has
no poles at the points of the discrete spectrum, we arrive at
\begin{equation}\label{4.10}
\frac{1}{2\pi\I}\oint_\sigma\,G_\pm(\lambda,x,y)d\lambda=\sum_{\lambda_k\in\sigma_d}\,\Res_{\lambda_k}\,
G_\pm^{\prime}(\lambda,x,y),\quad\pm y>\pm x.
\end{equation}
Moreover, the function $G^{\prime\prime}_\pm$ also does not contribute
to the left part of \eqref{4.10} since $G^{\prime\prime}_\pm(\lau,x,y)=
G^{\prime\prime}_\pm(\lal,x,y)$ for $\lambda\in\sigma^{(1)}_\mp$ and,
therefore, $\oint_{\sigma_\mp^{(1)}} G_\pm^{\prime\prime}(\lambda, x,y)
d\lambda=0$. In addition, $\oint_{\sigma_\pm} G_\pm^{\prime\prime}(\lambda,
x,y) d\lambda=0$ for $x\neq y$ by \eqref{1.14}.

Next we study  the contribution of the function $G_\pm^\prime$.
Properties \textbf{I, (b)} and \textbf{(c)} of Lemma \ref{lem2.3} imply
that
\begin{equation}\label{4.12}
|R_\pm(\lambda)|<1 \quad\mbox{for}\quad \lambda\in\inte(\sigma^{(2)}),\qquad
|R_\pm(\lambda)|=1\quad\mbox{for}\quad\lambda\in\sigma^{(1)}_\pm.
\end{equation}
Now we consecutively use \eqref{2.18}, \eqref{2.16},
\eqref{2.17}, \eqref{1.14}, \eqref{2.2}, \eqref{2.4} and again
\eqref{1.14}, to obtain
\begin{align} \nn
\frac{1}{2\pi\I} & \oint_{\sigma_\pm} G_\pm^\prime(\lambda,x,y) d\lambda=
\oint_{\sigma_\pm} T_\pm(\lambda) \phi_\mp(\lambda,x) \psi_\pm(\lambda,y)d\rho_\pm(\lambda)\\ \nn
& =  \oint_{\sigma_\pm} \left(R_\pm(\lambda)\phi_\pm(\lambda,x)
+\overline{\phi_\pm(\lambda,x)}\right) \psi_\pm(\lambda,y)d\rho_\pm(\lambda)\\ \nn
& =\oint_{\sigma_\pm} R_\pm(\lambda)\psi_\pm(\lambda,x)\psi_\pm(\lambda,y)d\rho_\pm(\lambda) +
\oint_{\sigma_\pm} \breve\psi_\pm(\lambda,x) \psi_\pm(\lambda,y)d\rho_\pm(\lambda)\\ \nn
& \quad \pm\int_x^{\pm\infty}d t\, K_\pm(x,t)\left(
\oint_{\sigma_\pm} R_\pm(\lambda)\psi_\pm(\lambda,t) \psi_\pm(\lambda,y)d\rho_\pm(\lambda)
+\delta(t-y)\right) \\ \label{4.13}
& =F_{r,\pm}(x,y) \pm \int_x^{\pm\infty} K_\pm(x,t)F_{r,\pm}(t,y)d t +
K_\pm(x,y),
\end{align}
where
\begin{equation}\label{4.14}
F_{r,\pm}(x,y)= \oint_{\sigma_\pm}R_\pm(\lambda) \psi_\pm(\lambda,x)\psi_\pm(\lambda,y)d\rho_\pm(\lambda).
\end{equation}
On the set $\sigma^{(1)}_\mp$ both the numerator and denominator of the
function $G_\pm^\prime$ have poles (resp., square root
singularities) at points of the set $\sigma^{(1)}_\mp\cap
(M_\pm\cup(\pa\sigma^{(1)}_+\cap\pa\sigma^{(1)}_-))$
(resp.\ $\sigma^{(1)}_\mp\cap (M_\mp\setminus(M_\mp\cap M_\pm))$), but
multiplying them, if necessary, by the functions \eqref{2.6}, we can avoid
singularities. Hence, without loss of generality, we can suppose
$\sigma^{(1)}_\mp\cap (M_{r_+}\cup M_{r_-})=\emptyset$.
Then, since $\psi_\pm(\lambda,x)\in\mathbb{R}$ as $\lambda\in\sigma_\mp^{(1)}$,
\begin{equation}\label{4.16}
\frac{1}{2\pi\I}\oint_{\sigma_\mp^{(1)}} G^\prime_\pm(\lambda,x,y)d\lambda
=\frac{1}{2\pi\I}\int_{\sigma^{(1),\mathrm{u}}_\mp}\psi_\pm(\lambda,y)\left(\frac{\overline{\phi_\mp(\lambda,x)}}{\overline{W(\lambda)}}
-\frac{\phi_\mp(\lambda,x)}{W(\lambda)}\right)d\lambda.
\end{equation}
Property \textbf{I, (b)} of Lemma \ref{lem2.3} and \eqref{2.16} imply
\[
\overline{\phi_\mp(\lambda,x)}=T_\mp(\lambda)\phi_\pm(\lambda,x) -
\frac{T_\mp(\lambda)}{\ov{T_\mp(\lambda)}}\phi_\mp(\lambda,x).
\]
Therefore,
\begin{align}\nn
\frac{\phi_\mp(\lambda,x)}{W(\lambda)} -\frac{\overline{\phi_\mp(\lambda,x)}}{\overline{W(\lambda)}}
& =
\phi_\mp(\lambda,x)\left(\frac{1}{W(\lambda)} +\frac{T_\mp(\lambda)}{\ov{T_\mp(\lambda)}\,\overline{W(\lambda)}}\right) -
\frac{T_\mp(\lambda)\phi_\pm(\lambda,x)}{\overline{W(\lambda)}}\\ \label{4.17}
& =\phi_\mp(\lambda,x)\frac{2\Re\left(T_\mp^{-1}(\lambda)\overline{W(\lambda)}\right)T_\mp(\lambda)}
{|W(\lambda)|^2}-
\frac{T_\mp(\lambda)\phi_\pm(\lambda,x)}{\overline{W(\lambda)}}.
\end{align}
But by \eqref{2.18}
\[
T_\mp^{-1}(\lambda)\overline{W(\lambda)} = |W(\lambda)|^2 g_\mp(\lambda) \in
\I\mathbb{R},\quad\mbox{for}\quad\lambda\in\sigma_\mp^{(1)},
\]
thus, the first summand in \eqref{4.17} vanishes. And using
$\overline{W}= (\overline{T_\mp} g_\mp)^{-1}$ we arrive at
\begin{equation}\label{4.18}
\frac{\overline{\phi_\mp(\lambda,x)}}{\overline{W(\lambda)}} - \frac{\phi_\mp(\lambda,x)}{W(\lambda)} =
|T_\mp(\lambda)|^2 g_\mp(\lambda) \phi_\pm(\lambda,x).
\end{equation}
Combining \eqref{4.18}, \eqref{1.12}, \eqref{4.16}, \eqref{2.2} and \eqref{4.13} we have
\begin{equation}\label{4.19}
\frac{1}{2\pi\I} \oint_{\sigma} G_\pm(\lambda,x,y)d\lambda= F_{c,\pm}(x,y) +
K_\pm(x,y)\pm\int_x^{\pm\infty}K_\pm(x,t)F_{c,\pm}(t,y)d t,
\end{equation}
where
\begin{align}\label{4.20}
F_{c,\pm}(x,y)& := F_{r,\pm}(x,y) + F_{h,\pm}(x,y),\\ \label{4.21}
F_{h,\pm}(x,y) &:= \int_{\sigma^{(1),\mathrm{u}}_\mp} \psi_\pm(\lambda,x)
\psi_\pm(\lambda,y) |T_\mp(\lambda)|^2d\rho_\mp(\lambda).
\end{align}
To derive the part of the GLM equation kernel, that correspond to the discrete
spectrum (for function $F_{c,\pm}$ index $c$ means the part,
corresponding to the continuous spectrum), we  apply \eqref{2.12},
\eqref{2.10}, \eqref{D.1}, \eqref{D.3} and \eqref{2.2} to the right
hand side of \eqref{4.10}. Then,
\begin{align} \nn
&
\sum_{\lambda_k\in\sigma_d}\,\Res_{\lambda_k}\, G_\pm^{\prime}(\lambda,x,y)
= -\sum_{\lambda_k\in\sigma_d}\,\Res_{\lambda_k}
\frac{\tilde\phi_\mp(\lambda,x)\tilde\psi_\pm(\lambda,y)}{\tilde W(\lambda)}\\ \nn
&
= -\sum_{\lambda_k\in\sigma_d} \frac{\tilde\phi_\pm(\lambda_k,x)\tilde\psi_\pm(\lambda_k,y)}{\tilde
W^\prime(\lambda_k)c_{\pm,k}} =
-\sum_{\lambda_k\in\sigma_d}\,(\gamma_k^\pm)^2\tilde\phi_\pm(\lambda_k,x)\tilde\psi_\pm(\lambda_k,y)\\ \label{4.23}
&
= - F_{d,\pm}(x,y)\mp\int_x^{\pm\infty}K_\pm(x,t)F_{d,\pm}(t,y) dt,
\end{align}
where
\begin{equation}\label{4.24}
F_{d,\pm}(x,y):=\sum_{\lambda_k\in\sigma_d}(\gamma_k^\pm)^2\tilde\psi_\pm(\lambda_k,x)\tilde\psi_\pm(\lambda_k,y).
\end{equation}
Combining \eqref{4.10}, \eqref{4.19}, and \eqref{4.23} we finally obtain \eqref{4.2}.
\end{proof}

As is shown in Appendix~\ref{secTOP}, the kernel $F_\pm(x,y)$ of the GLM equation satisfies

\begin{lemma}\label{lem4.1}
The kernel of the GLM equation possess the following properties:
\begin{enumerate}[\bf I.]
\addtocounter{enumi}{3}
\item
The function $F_\pm(x,y)$ is continuously differentiable with respect to both variables and
there exists real-valued function $q_\pm(x)$, $x\in\mathbb{R}$, with
\[
\pm\int_a^{\pm\infty}(1+x^2) |q_\pm(x)| dx < \infty,  \quad\mbox{for all } a\in\mathbb{R},
\]
such that
\begin{align}\label{4.3}
& |F_\pm(x,y)|\leq C_\pm(x)
Q_\pm\left(x+y\right),\\ \label{4.4}
& \left|\frac{\pa}{\pa x} F_\pm (x,y)\right|\leq
C_\pm(x)\left(\left|q_\pm\left(\frac{x+y}{2}\right)\right|+Q_\pm(x+y)\right),\\ \label{4.5}
& \pm\int_a^{\pm\infty}\,\left|\frac{d}{dx} F_\pm (x,x)\right|(1+ x^2)\,dx<\infty,
\end{align}
where
\begin{equation}\label{4.31}
Q_\pm(x):=\pm\int_{\frac{x}{2}}^{\pm\infty}\left|q_\pm(t)\right|d t,
\end{equation}
and $C_\pm(x)>0$ is a continuous function, which decreases monotonically as $x\to\pm\infty$.
\end{enumerate}
\end{lemma}

In summary, we have obtained the following necessary conditions
for the scattering data:

\begin{theorem}[necessary conditions for the scattering data]\label{theor1}
The scattering data
\begin{align}\nn
{\mathcal S} = \Big\{ & R_+(\lambda),\,T_+(\lambda),\, \lambda\in\sigma_+^{\mathrm{u,l}}; \,
R_-(\lambda),\,T_-(\lambda),\, \lambda\in\sigma_-^{\mathrm{u,l}};\\\label{4.6}
& \lambda_1,\dots,\lambda_p\in\mathbb{R}\setminus (\sigma_+\cup\sigma_-),\,
\gamma_1^\pm,\dots,\gamma_p^\pm\in\mathbb{R}_+\Big\}
\end{align}
possess the properties \emph{\textbf{I-III}} listed in Lemma~\ref{lem2.3}. The
functions $F_\pm(x,y)$, defined in \eqref{4.2}, possess property
\emph{\textbf{IV}} from Lemma~\ref{lem4.1}.
\end{theorem}

In fact, the conditions on the scattering data, given in this theorem are both
necessary and sufficient for the solution of the scattering problem in the class
\eqref{0.1}--\eqref{0.2}. The sufficiency of these conditions together with the
algorithm for the solution of the inverse problem will be discussed in the next section.

As a consequence of the GLM equation and its unique solvability (see Lemma~\ref{lem5.4}) and
also formula \eqref{A.19} we note

\begin{corollary}\label{cor-theor1}
Suppose $q(x)$ satisfies \eqref{0.2}, then $q(x)$ is uniquely determined by
one of the sets of its ``partial" scattering data ${\mathcal S}_+$ or ${\mathcal S}_-$, where
\begin{align}\nn
{\mathcal S}_\pm = \Big\{ & R_\pm(\lambda),\, \lambda\in\sigma_\pm^{\mathrm{u}}; \,
|T_\mp(\lambda)|^2,\, \lambda\in\sigma_\mp^{(1),\mathrm{u}};\\
& \lambda_1,\dots,\lambda_p\in\mathbb{R}\setminus (\sigma_+\cup\sigma_-),\,
\gamma_1^\pm,\dots,\gamma_p^\pm\in\mathbb{R}_+\Big\}.
\end{align}
\end{corollary}

The question about the characterization of such sets (necessary and sufficient conditions) for potentials from the class \eqref{0.2} is rather complicated and is still open.

\section{The inverse scattering problem}           \label{secISP}

Let $H_\pm$ be two one-dimensional finite-gap Schr\"odinger
operators associated with potentials $p_\pm(x)$ as introduced in
Section~\ref{secBgOp}. Let $\mathcal{S}$ be given scattering data
with corresponding kernels $F_\pm(x,y)$ satisfying the necessary
conditions from Theorem~\ref{theor1}.

We begin by showing that, given $F_\pm(x,y)$, the GLM equations \eqref{ME} can be
solved for $K_\pm(x,y)$.

\begin{lemma}\label{lem5.4}
Under condition \emph{\textbf{IV}}, \eqref{4.3}, the GLM equations \eqref{ME}
have unique real-valued solutions $K_\pm(x,\,\cdot\,)\in L^1(x,\pm\infty)$
satisfying the estimates
\begin{equation}\label{5.100}
|K_\pm(x,y)|\leq C_\pm(x)Q_\pm(x+y),\quad\pm y>\pm x.
\end{equation}
Here the functions $Q_\pm(x)$ are the same, as in
\eqref{4.31}, and $C_\pm(x)$ are functions of the same type, as
in Lemma~\ref{lem4.1} (i.e.\ positive, continuous and decreasing as
$x\to\pm\infty$).

Moreover, under condition \emph{\textbf{IV}}, \eqref{4.4},
the functions $K_\pm(x,y)$ are differentiable with respect to each
variable and satisfy the estimate \eqref{A.22}, where the
functions $q_\pm(x)$ are the same as in \eqref{4.4} and the
functions $C_\pm(x)$ are of the same type  as in \eqref{5.100}.
Besides,
\begin{equation}\label{5.101}
\pm\int_a^{\pm\infty}\,(1+x^2)\left|\frac{d}{dx}\,K_\pm(x,x)\right| dx<\infty,
\qquad \forall a\in\mathbb{R}.
\end{equation}
\end{lemma}

\begin{proof}
The solvability of \eqref{ME} under condition \eqref{4.3} together with the
estimate \eqref{5.100} follows from considerations completely analogous to
those ones used in the proof of Lemma~\ref{lemA.3} (see Remark~\ref{rem6}).
To prove uniqueness, first note that the GLM equations are generated by
compact operators. Thus it is sufficient to prove, that the equation
\begin{equation}\label{5.102}
f(x)\pm\int_x^{\pm\infty} F_\pm(x,y) f(y)d y =0
\end{equation}
has only the trivial solution in the space $L^1(x,\pm\infty)$.
The proof is similar for the ``$+$" and ``$-$"
cases, hence we give it only for the ``$+$" case. Let $f(y)$, $y>x$, be a
nontrivial solution of \eqref{5.102} and set $f(y)=0$ for $y\leq x$.
Since $F_+(x,y)$ is real-valued, we can assume $f(y)$ is real-valued.
Abbreviate by
\begin{equation}\label{5.103}
\widehat f(\lambda) =\int_\mathbb{R} \psi_+(\lambda,y) f(y)d y
\end{equation}
the generalized Fourier transform, generated by the spectral decomposition \eqref{1.14}
(cf.~\cite{T}). Recall that $\widehat f(\lambda)\in L_{\loc}^1 (\siu_+\cup\sil_+)$ and $\widehat f(\lambda)=O(1)$
as $\lambda\to+\infty$.

Multiplying \eqref{5.102} by $f(x)$, integrating over $\mathbb{R}$, and
applying \eqref{5.103} and \eqref{4.2} we have
\begin{align}\nn
& 2 \int_{\sigma_+^u}|\widehat f(\lambda)|^2d\rho_+(\lambda)
+2\Re \int_{\sigma_+^u}R_+(\lambda) \widehat f(\lambda)^2d\rho_+(\lambda)\\ \label{5.104}
& \quad +\int_{\sigma_-^{(1),u}} \widehat f(\lambda)^2 |T_-(\lambda)|^2d\rho_-(\lambda) +
\sum_{k=1}^p (\gamma ^+_k)^2 \left(\int_\mathbb{R} \tilde\psi_+(\lambda_k,y) f(y)d y \right)^2=0.
\end{align}
The last two summands in \eqref{5.104} are nonnegative since $\widehat f(\lambda)\in\mathbb{R}$
for $\lambda\in\sigma_-^{(1)}$ and $\tilde\psi_+(\lambda_k,x)\in\mathbb{R}$. Ignoring the last one
and proceeding as in \cite[Lemma~3.5.3]{M} we obtain
\begin{equation}\label{5.105}
2\int_{\sigma^{(2),\mathrm{u}}} (1-|R_+(\lambda)|)|\widehat f(\lambda)|^2d\rho_+(\lambda) +
\int_{\sigma_-^{(1),\mathrm{u}}} \widehat f(\lambda)^2 |T_-(\lambda)|^2d\rho_-(\lambda) \leq 0.
\end{equation}
Here we used that
\[
\int_{\sigma_+^{(1),\mathrm{u}}}(1-|R_+(\lambda)|) |\widehat f(\lambda)|^2d\rho_+(\lambda) = 0
\]
by condition \textbf{I, (b)}. Now, since $|R_+(\lambda)|<1$, $\rho_+(\lambda)>0$ for $\lambda\in \inte(\sigma^{(2)})$
and $\rho_-(\lambda)>0$ for $\lambda\in\inte(\sigma_-^{(1)})$, we conclude that
\[
\widehat f(\lambda)=0 \quad \text{for}\quad \lambda\in\sigma^{(2)}\cup \sigma_-^{(1)}= \sigma_-.
\]
The function $\widehat f(z)$ can be defined by formula \eqref{5.103} as a
meromorphic function on $\mathbb{C}\setminus\sigma_+$. By our analysis it is even
meromorphic on $\mathbb{C}\setminus\sigma^{(1)}_+$ and vanishes on $\sigma_-$,
thus $\widehat f(z)$ is equal to zero and hence also $f(x)$.

The estimate \eqref{5.101} follows by literally repeating the proof of
Lemma~\ref{lemA.3}.
\end{proof}

Next, define two functions
\begin{equation}\label{5.1}
\tilde q_\pm(x) =\mp\frac{d^2}{dx^2}K_\pm(x,x) +
p_\pm(x),\quad x\in\mathbb{R}
\end{equation}
and note that estimate \eqref{5.101} implies
\begin{equation}
\pm \int_a^{\pm \infty}|\tilde q_\pm(x) - p_\pm(x)|
(1+x^2) d x <\infty ,\quad a\in\mathbb{R}. \label{5.2}
\end{equation}

\begin{lemma}\label{lem5.7}
The functions $\phi_\pm(z,x)$, defined by
\begin{equation}\label{Dop}
\phi_\pm(z,x) =\psi_\pm(z,x)\pm\int_{x}^{\pm\infty}\,K_\pm(x,y)\psi_\pm(z,y)\,d y,
\end{equation}
solve the equations
\begin{equation}\label{5.3}
\left(-\frac{d^2}{dx^2}
+\tilde q_\pm(x)\right)\,\phi_\pm(z,x) = z\phi_\pm(z,x),
\end{equation}
where $\tilde q_\pm(x)$ are defined by \eqref{5.1}.
\end{lemma}

\begin{proof}
Consider the two operators\footnote{We don't know $\tilde H_\pm$ is limit point
at $\mp\infty$ yet, but this will not be used.}
\[
\tilde H_\pm = -\frac{d^2}{d x^2} +\tilde q_\pm(x),\quad x\in\mathbb{R}.
\]
On the corresponding half-axes  the potentials $\tilde q_\pm(x)$ are
asymptotically close to our background potentials $p_\pm(x)$.
Define the integral operators
\[
\left ({\mathcal K}_\pm f\right)(x) =\pm\int_x^{\pm\infty} K_\pm(x,y) f(y)d y.
\]
Under the assumption, that the kernel $F_\pm(x,y)$ of the GLM equation is
twice continuously differentiable, we infer from \eqref{ME} that the function
$K_\pm(x,y)$ is also twice differentiable (see the proof of Lemma~\ref{lemA.3}).
Moreover one can prove, literally following \cite{F2}, that the identity
\[
\tilde H_\pm{\mathcal K}_\pm={\mathcal K}_\pm\,H_\pm,
\]
is valid. This identity implies \eqref{5.3}.  To obtain equality \eqref{5.3} without
assumption of existence of the second derivatives, one can literally
follow the proof of  \cite[Theorem~3.3.1]{M}.
\end{proof}

The remaining problem is to show $\tilde q_+(x)\equiv \tilde q_-(x)$
under conditions \textbf{II} and \textbf{III} on the scattering data $\mathcal{S}$.

\begin{theorem}[uniqueness of the reconstructed potential]\label{theor2}
Let the scattering data ${\mathcal S}$, defined as in \eqref{4.6},
satisfy conditions \emph{\textbf{I, (a)--(d)}, \textbf{II}, \textbf{III, (a)}} and \emph{\textbf{IV}}. Then each
of the GLM equations \eqref{ME} has a unique solution
$K_\pm(x,y)$, satisfying the estimate \eqref{5.101}.
The functions $\tilde
q_\pm(x)$, defined by \eqref{5.1}, satisfy \eqref{5.2}.

Under additional conditions \emph{\textbf{III, (b)}} and \emph{\textbf{I, (e)}}, these two functions
coincide, $\tilde q_-(x)\equiv \tilde q_+(x)=:q(x)$, and the data
${\mathcal S}$ are the scattering data for the Schr\"odinger
operator with potential $q(x)$.
\end{theorem}

\begin{proof}
To prove the uniqueness of the reconstructed potential we follow the method
proposed in \cite{M}. Recall that, according to \cite{L,T}, we have the
inversion formula for the generalized Fourier transform, generated
by the spectral decomposition \eqref{1.14} and applied to the function
$f(\lambda)\in L_{\loc}^1 (\sipmu\cup\sipml)$, $f(\lambda)=O(1)$, $\lambda\to+\infty$:
\begin{align}\nn
\check f(y) &=\oint_{\sigma_\pm} f(\lambda) \psi_\pm(\lambda,y)d\rho_\pm(\lambda),\\ \label{5.4}
f(\lambda) &=\int_\mathbb{R}\check f(y)\overline{\psi_\pm(\lambda,y)}d y.
\end{align}
Split the kernel of the GLM equation \eqref{4.2} according to $F_\pm(x,y)=F_{r,\pm}(x,y) +
F_{h,\pm}(x,y) +F_{d,\pm}(x,y)$ (cf.\ \eqref{4.14}, \eqref{4.21}, \eqref{4.24}).

We begin by considering the following part of the GLM equation
\begin{equation}\label{5.61}
G_\pm(x,y):=F_{r,\pm}(x,y)
\pm\int_x^{\pm\infty}K_\pm(x,t)F_{r,\pm}(t,y)d t,
\end{equation}
where
$K_\pm(x,y)$ are the solutions of GLM equations. By condition
\textbf{I, (b)--(c)} we have $|R_\pm(\lambda)|\leq 1$ and properties (iii)
and (i) of Lemma~\ref{lem1.1} imply, that we can take $f(\lambda)=R_\pm(\lambda)
 \psi_\pm(\lambda,x)$ in \eqref{5.4}. Using \eqref{4.14} we obtain
\begin{equation}\label{5.6}
\int_\mathbb{R}
F_{r,\pm}(x,y) \overline{\psi_\pm(\lambda,y)}d y =R_\pm(\lambda)\psi_\pm(\lambda,x).
\end{equation}
and \eqref{2.2} consequently implies
\begin{equation}\label{5.7}
\int_\mathbb{R}G_\pm(x,y) \overline{\psi_\pm(\lambda,y)}d y =R_\pm(\lambda)\phi_\pm(\lambda,x),
\quad\lambda\in\sipmul.
\end{equation}
On the other hand, invoking the GLM equations we have for $\pm y>\pm x$,
\begin{align*}
G_\pm(x,y)
&=  -K_\pm(x,y)
  -F_{h,\pm}(x,y)  -  F_{d,\pm}(x,y)\\
&\quad
\mp \int_{\sigma_\mp^{(1),\mathrm{u}}}d\rho_\mp(\lambda) |T_\mp(\lambda)|^2
\psi_\pm(\lambda,y) \int_x^{\pm\infty}K_\pm(x,t) \psi_\pm(\lambda,t)d t\\
&\quad
\mp
\sum_{k=1}^p (\gamma_k^\pm)^2\tilde\psi_\pm(\lambda_k,y) \int_x^{\pm\infty}K_\pm(x,t)
\tilde\psi_\pm(\lambda_k,t)d t.
\end{align*}
Again using \eqref{2.2} this gives
\begin{align} \nn
G_\pm(x,y)
&= -K_\pm(x,y) - \int_{\sigma_\mp^{(1),\mathrm{u}}}
|T_\mp(\lambda)|^2 \psi_\pm(\lambda,y) \phi_\pm(\lambda,x)d\rho_\mp(\lambda) \\ \label{5.9}
&\quad - \sum_{k=1}^p (\gamma_k^\pm)^2\tilde\psi_\pm(\lambda_k,y) \tilde\phi_\pm(\lambda_k,x).
\end{align}
Now we use this formula to evaluate
\begin{align*}
\int_{\mathbb{R}} G_\pm(x,y) \breve\psi_\pm(\lambda,y)d y= &
\mp\int_x^{\mp\infty}G_\pm(x,y) \breve\psi_\pm(\lambda,y)d y
\mp\int_x^{\pm\infty}K_\pm(x,y)\breve\psi_\pm(\lambda,y)d y\\
& \pm\int_{\sigma_\mp^{(1),\mathrm{u}}} |T_\mp(\xi)|^2
\phi_\pm(\xi,x)\,W(\psi_\pm(\xi,x),\breve\psi_\pm(\lambda,x))\frac{d\rho_\mp(\xi)}{\xi - \lambda}\\
& \pm
\sum_{k=1}^p\,(\gamma_k^\pm)^2\tilde\phi_\pm(\lambda_k,x)\,\frac{W(\tilde\psi_\pm(\lambda_k,x),
\breve\psi_\pm(\lambda,x))}{\lambda_k - \lambda}.
\end{align*}
Here we have used
\[
\left(\xi -
\lambda\right)\int_x^{\pm\infty}\breve\psi_\pm(\lambda,y) \psi_\pm(\xi,y)d y=
W(\breve\psi_\pm(\lambda,x),\psi_\pm(\xi,x)),
\]
which follows from Green's formula and the fact, that $\tilde\psi_\pm(\xi,y)\to 0$,
$\tilde\psi_\pm^\prime(\xi,y)\to 0$ for $\xi\notin\sigma_\pm$ as $y\to\pm\infty$.

Combining this formula with \eqref{5.7} and using \eqref{2.4} we infer the relation
\begin{equation}\label{5.10}
R_\pm(\lambda)\,\phi_\pm(\lambda,x) +
\overline{\phi_\pm(\lambda,x)} =T_\pm(\lambda)\theta_\mp(\lambda,x),
\quad\lambda\in\sipmul,
\end{equation}
where
\begin{align}\label{5.11}
\theta_\mp(\lambda,x) &:=\frac{1}{T_\pm(\lambda)}\left(
\breve\psi_\pm(\lambda,x) \mp\int_x^{\mp\infty}G_\pm(x,y)
\breve\psi_\pm(\lambda,y)d y \right.\\ \nn &\qquad
-\int_{\sigma_\mp^{(1),\mathrm{u}}} |T_\mp(\xi)|^2
\phi_\pm(\xi,x)W(\psi_\pm(\xi,x),\breve\psi_\pm(\lambda,x))
\frac{d\rho_\mp(\xi)}{\xi - \lambda}\\ \nn &\qquad
\left.\pm\sum_{k=1}^p (\gamma_k^\pm)^2 \tilde\phi_\pm(\lambda_k,x)
\frac{W(\tilde\psi_\pm(\lambda_k,x),\breve\psi_\pm(\lambda,x))}{\lambda_k
- \lambda}\right).
\end{align}
It turns out that, in spite of the fact that $\theta_\mp(\lambda,x)$ is defined via the
background solutions corresponding to the opposite half-axis $\mathbb{R}_\pm$, it
shares a series of properties with $\phi_\mp(\lambda,x)$. Namely, we prove

\begin{lemma}\label{lem5.1}
Let $\theta_\mp(z,x)$ be defined by formula \eqref{5.11} on the set $\sipmul$.
\begin{enumerate}[\rm(i)]
\item
The function $\tilde\theta_\mp(z,x)= \delta_\mp(z) \theta_\mp(z,x)$
admits an analytical extension to the domain $\mathbb{C}\setminus\sigma$.
\item
The function $\tilde\theta_\mp(z,x)$ is continuous up to the
boundary $\siul$ except possibly at the points
$\pa\sigma_+\cup\pa\sigma_-$. Furthermore,
\begin{equation}\label{5.30}
\theta_\mp(\lau,x)=
\begin{cases}
\theta_\mp(\lal,x)\in\mathbb{R},&\text{for }
\lambda\in\mathbb{R}\setminus\sigma_\mp,\\[1mm]
\overline{\theta_\mp(\lal,x)},&\text{for }\lambda\in\inte(\sigma_\mp).
\end{cases}
\end{equation}
\item
For large $z$ the function $\theta_\mp(z,x)$ has the following asymptotic
behavior
\[
\theta_\mp(z,x)=\E^{\mp \I\sqrt{z}\,x}\bigl(1+O(z^{-1/2})\bigr),\ \ z\to\infty.
\]
\item
$W\bigl(\theta_\mp(z,x),\phi_\pm(z,x)\bigr)=\pm W(z)$, where
$W(z)$ is defined by formula \eqref{2.18}.
\end{enumerate}
\end{lemma}

\begin{remark}\label{rem5}
Note that we did not establish the connection between the function
$W(z)$ and the function $W\bigl(\phi_+(z,x),
\phi_-(z,x)\bigr)$, which can depend on $x$, because $\phi_+$
and $\phi_-$ are the solutions of  Schr\"odinger equations
corresponding to possibly different potentials $\tilde q_+$ and
$\tilde q_-$.
\end{remark}

\begin{proof}
To show (i) we will show that $\tilde\theta_\mp(z,x)$
has an analytic extension to $\mathbb{C}\setminus\sigma$.
We will study each term in \eqref{5.11} separately.

For the first one,
\begin{equation}\label{5.141}
\zeta_\mp(z,x):=\frac{\breve\psi_\pm(z,x)}{T_\pm(z)},
\end{equation}
it is immediate by \eqref{2.19} that $\tilde\zeta_\mp(z,x)= \delta_\mp(z) \zeta_\mp(z,x)$ has the
required property. This also covers the second term since
$G_\pm(x,\,\cdot\,)\in L^1(\mathbb{R})$ is real-valued.

Now we discuss the properties of the Cauchy-type integral in the
representation \eqref{5.11}. Multiplying it by $T_\pm^{-1}(z)$, we
represent the third summand in \eqref{5.11} as
\begin{equation}\label{5.123}
H_\mp(z,x):=\mp \frac{1}{2\pi\I} \int_{\sigma_\mp^{(1),\mathrm{u}}}
h_\mp(z,\xi,x)\frac{d\xi}{\xi-z},
\end{equation}
where the integrand, due to \eqref{2.10}, \eqref{1.12}, and \eqref{2.18},
has the representation
\begin{align}\nn
h_\mp(z,\xi,x)
&= \frac{\delta_\mp(\xi)^2}{g_\mp(\xi)
|\tilde W(\xi)|^2} \tilde\phi_\pm(\xi,x)W(\tilde\psi_\pm(\xi,x),\zeta_\mp(z,x))\\ \label{5.13}
& = \frac{|\hat\delta_\mp(\xi)|^2}{g_\mp(\xi)
|\hat W(\xi)|^2} \frac{|\hat\delta_\pm(\xi)|^2}{\hat\delta_\pm(\xi)^2}
\hat\phi_\pm(\xi,x)W(\hat\psi_\pm(\xi,x), \zeta_\mp(z,x)).
\end{align}
By property \textbf{II, (b)} the function $\hat W(\xi)$ has no zeros in the interior of
$\sigma_\mp^{(1),\mathrm{u}}$. Thus, for $z\notin \sigma_\mp^{(1)}$, the function $h_\mp(z,.,x)$
is bounded in the interior of $\sigma_\mp^{(1)}$ and the only possible singularities can arise
at  the boundary. We claim
\begin{equation}\label{5.131}
h_\mp(z,\xi,x)=
\begin{cases}
O(\sqrt{\xi -E}) &\text{for } E\notin\sigma_v,\\
O\left(\frac{1}{\sqrt{\xi -E}}\right)&\text{for } E\in\sigma_v,
\end{cases}
\qquad E\in\pa\sigma_\mp^{(1)},\; z\neq E.
\end{equation}
This follows from $\frac{|\hat\delta_\mp(\xi)|^2}{g_\mp(\xi)}=O(\sqrt{\xi -E})$
together with $\hat W(\xi) = O(1)$ if $E\notin\sigma_v$ and $\hat W(\xi) = C(E)(\sqrt{\xi -E})(1+o(1))$, $C(E)\neq 0$ by \textbf{II, (b)} if $E\in\sigma_v$.

So $h_\mp$ is integrable and the third summand of \eqref{5.11} also inherits the properties
of $\zeta_\mp(z,x)$.

Finally, let us consider the last summand in \eqref{5.11}. It again
inherits everything from $\tilde\zeta_\mp(z,x)$ except for possible
additional poles at the points $\lambda_k$. However, these are canceled
by the fact that the function $\tilde W(z)$ vanishes for $z=\lambda_k$.

(ii).
Next we look at the boundary values. The only nontrivial term
is of course the Cauchy-type integral \eqref{5.123} as
$z\to\lambda\in\inte(\sigma_\mp^{(1)})$. First of all observe that by
\eqref{1.62} and \eqref{2.18} we have
\[
\frac{W(\tilde\psi_\pm(\lambda,x),\breve\psi_\pm(z,x))}{T_\pm(z)} \to
\pm (\delta_\pm W)(\lambda),
\]
where the function $\delta_\pm W$ is bounded and non zero for $\lambda\in\inte(\sigma_\mp^{(1)})$ by
\textbf{II, (a)}. Hence the Plemelj formula applied to \eqref{5.123} gives
\[
H_\mp(\lambda,x) = \pm
\frac{\tilde\phi_\pm(\lambda,x)}{2\delta_\pm(\lambda) g_\mp(\lambda) \overline{W(\lambda)}}
\mp \dashint_{\sigma_\mp^{(1), \mathrm{u}}}\frac{h_\mp(\lambda,\xi,x)}{\xi-\lambda}d\xi,
\quad \lambda\in\inte(\sigma_\mp^{(1),\mathrm{u}}),
\]
where both terms are finite. Here $\dashint$ denotes the principal value integral.

Hence the boundary values away from $\pa\sigma_+\cup\pa\sigma_-$ exist and we have
\begin{equation} \label{5.14}
\theta_\mp(\lau,x)= \overline{\theta_\mp(\lal,x)}, \quad\lambda\in\sigma_+\cup\sigma_-.
\end{equation}
Moreover, by property \textbf{I, (b)} we have
\begin{equation}\label{5.16}
\theta_\mp=T_\pm^{-1}\left(R_\pm\phi_\pm +\overline{\phi_\pm}\right) =
\frac{\phi_\pm}{\ov{T_\pm}} +\frac{\overline{\phi_\pm}}{T_\pm}\in\mathbb{R}
\quad\mbox{for}\quad  \lambda\in\sigma_\pm^{(1)},
\end{equation}
from which
\begin{equation}\label{5.15}
\theta_\mp(\lau,x)= \theta_\mp(\lal,x) ,\quad \lambda\in\sigma_\pm^{(1)},
\end{equation}
follows. Combining \eqref{5.14} and \eqref{5.15} we have \eqref{5.30}.

(iii). For $|z|\to\infty$ due to properties (iii) of Lemma \ref{lem1.1} and \textbf{I, (d)} we have
\begin{equation}\label{5.18}
\zeta_\mp(z,x)= \E^{\mp \I\sqrt{z} x}\left(1+O(z^{-1/2})\right).
\end{equation}
Since the last two terms in \eqref{5.11} are $O(z^{-1})$ we obtain
\[
\theta_\mp(z,x)= \E^{\mp \I\sqrt{z} x} \left(1 + \int_0^\infty G_\pm(x,x\mp t)
\E^{\I\sqrt{z} t}d t +O(z^{-1/2})\right)
\]
which implies (iii) since $G_\pm(x,y)$ is differentiable with respect to $y$.

(iv). From \eqref{5.10} (invoking \eqref{2.18}) we obtain
\begin{equation}
W\bigl(\theta_\mp(z,x), \phi_\pm(z,x)\bigr)=\pm W(z)
\end{equation}
for $z\in \sigma_\pm$. Hence equality holds for all $z\in\mathbb{C}$ by analytical continuation.
\end{proof}

\begin{corollary}\label{col7}
The function $\tilde\theta_\mp(z,x)$ admits an analytical extension to the domain
$\mathbb{C}\setminus\sigma_\mp$.
\end{corollary}

\begin{proof}
Property (i) holds for $z\in\mathbb{C}\setminus \sigma$. Relation
\eqref{5.30} implies that $\tilde\theta_\mp$ has no jump across
$z\in\inte(\sigma_\pm^{(1)})$. To finish the proof we need to show that
the possible remaining singularities at $E\in\pa\sigma_\pm^{(1)}\cap\pa\sigma$
are removable. This follows from (cf.\ \eqref{2.19}) 
\begin{equation}\label{5.142}
\hat\zeta_\mp(z,x)= \frac{\hat W(z)}{2\sqrt{\prod_{j=0}^{2r_\pm} (z-E_j^\pm)}}
\breve\delta_\pm(z) \breve\psi_\pm(z,x)
\end{equation}
which shows $\tilde\zeta_\mp(z,x)=O((z-E)^{-1/2})$ and hence
$\tilde\theta_\mp(z,x)= O((z-E)^{-1/2})$ for $E\in\sigma_\pm^{(1)}\cap\pa\sigma$.

However, let us emphasize at this point that the behavior of
$\theta_\pm(z,x)$ at the remaining edges is a more subtle question
to be discussed later.
\end{proof}

Eliminating $\overline{\phi_\pm}$ from
\[
\begin{cases}
\overline{ R_\pm(\lambda)}\,\overline{\phi_\pm(\lambda,x)} + \phi_\pm(\lambda,x)=\overline{\theta_\mp(\lambda,x)}\,
\overline{T_\pm(\lambda)}&\\[2mm]
R_\pm(\lambda)\,\phi_\pm(\lambda,x) + \overline{\phi_\pm(\lambda,x)}=\theta_\mp(\lambda,x)\,T_\pm(\lambda)&
\end{cases}
\]
we obtain
\[
\phi_\pm(\lambda,x)\left(1 - |R_\pm(\lambda)|^2\right) =\overline{\theta_\mp(\lambda,x)}\,
\overline{T_\pm(\lambda)} -
\overline{R_\pm(\lambda)}\,\theta_\mp(\lambda,x)\,T_\pm(\lambda).
\]
Next, using \textbf{I, (c)}, \textbf{II} and the consistency condition \textbf{I, (e)} then shows
\begin{align*}
T_\mp(\lambda)\phi_\pm(\lambda,x) &=\overline{\theta_\mp(\lambda,x)} - \frac{\overline{R_\pm(\lambda)}
T_\pm(\lambda)}{\ov{T_\pm(\lambda)}} \theta_\mp(\lambda,x)\\
&=\overline{\theta_\mp(\lambda,x)} + R_\mp(\lambda) \theta_\mp(\lambda,x),
\quad\lambda\in\sigma^{(2)}.
\end{align*}
This equation together with \eqref{5.10} gives us a system from
which we can eliminate the reflection coefficients $R_\pm$. We obtain
\begin{equation}\label{5.20}
T_\pm(\lambda) (\phi_\pm(\lambda)\phi_\mp(\lambda) -
\theta_\pm(\lambda)\theta_\mp(\lambda))=
\phi_\pm(\lambda)\overline{\theta_\pm(\lambda)} -
\overline{\phi_\pm(\lambda)}\theta_\pm(\lambda), \quad
\lambda\in\sigma^{(2),\mathrm{u,l}}.
\end{equation}
Now introduce the function
\begin{equation}\label{5.211}
G(z):= G(z,x) = \frac{\phi_+(z,x)\phi_-(z,x) - \theta_+(z,x)\theta_-(z,x)}{W(z)}
\end{equation}
which is well defined in the domain $z\in \mathbb{C}\setminus
\left(\sigma\cup\sigma_d\cup M_+\cup M_- \right)$.

From \eqref{5.20} and \eqref{2.18} we see, that
\begin{equation}\label{5.22}
G(\lambda) = -\left(\phi_\pm(\lambda)\overline{\theta_\pm(\lambda)}
- \overline{\phi_\pm(\lambda)}\theta_\pm(\lambda)\right)
g_\pm(\lambda), \qquad\lambda\in\sigma^{(2),\mathrm{u,l}}.
\end{equation}
So we need to study the properties of $G(z,x)$ as a function of $z$ (regarding $x$ as
a fixed parameter). Our aim  is to prove that $G(z,x)=0$. This will follow from

\begin{lemma}\label{lem5.2}
The function $G(z,x)$ possess the following properties.
\begin{enumerate}[\rm(i)]
\item
\begin{equation}\label{5.264}
G(\lau,x) =  G(\lal,x)\in\mathbb{R} \ \mbox{for}\
\lambda\in\mathbb{R}\setminus (\pa\sigma_-\cup\pa\sigma_+\cup\sigma_d).
\end{equation}
\item
It has removable singularities at the points $\pa\sigma_-\cup\pa\sigma_+\cup\sigma_d$,
where $\sigma_d:=\{\lambda_1,\dots,\lambda_p\}$.
\end{enumerate}
\end{lemma}

\begin{proof}
(i). We can rewrite $G(z,x)$ as
\begin{equation}\label{5.25}
G(z,x)= \frac{\tilde\phi_+(z,x)\tilde\phi_-(z,x) -
\tilde\theta_+(z,x)\tilde\theta_-(z,x)}{\tilde W(z)},
\end{equation}
where $\tilde\theta_\pm(z,x)=\delta_\pm(z)\theta_\pm(z,x)$ as usual.
The numerator is bounded near the points under consideration, and the
denominator does not vanish there.
Thus $G(z,x)$ has no singularities at the points $(M_+\cup M_-) \setminus \sigma_d$.

Furthermore, by Lemma~\ref{lem5.1}, \textbf{II, (a)}, and Lemma~\ref{lem2.1}, (i)
we know that $G(z,x)$ has continuous limiting values on the sets
$\sigma_-$ and $\sigma_+$, except possibly at the edges, satisfying
\[
G(\lau,x) =\overline{G(\lal,x)},\quad \lambda\in\sigma_+ \cup \sigma_-.
\]
Hence, if we can show that these limits are real, they will be equal and
$G(z,x)$ will extend to a meromorphic function on $\mathbb{C}$, that is, (i) holds.

First of all observe that from \eqref{5.30}, \eqref{5.22}, \eqref{1.8}, and Lemma~\ref{lem2.1} (i),
it follows, that
\begin{equation}\label{5.24}
G(\lau,x) = G(\lal,x)\in\mathbb{R}, \quad \lambda\in\inte(\sigma^{(2)}).
\end{equation}
Thus, it remains to prove
\begin{equation}\label{5.26}
G(\lau,x) =  G(\lal,x)\in\mathbb{R}
\quad \mbox{for}\quad \lambda\in\inte(\sigma_-^{(1)})\cup\inte(\sigma_+^{(1)}).
\end{equation}
So let us show that there is no jump on the set
$\inte(\sigma_-^{(1)})\cup\inte(\sigma_+^{(1)})$. Abbreviate
\begin{equation}\label{5.28}
\left[G\right] :=G(\lambda)
-\overline{G(\lambda)}=\left[\frac{\phi_+\phi_-}{W}\right] -
\left[\frac{\theta_+\theta_-}{W}\right], \qquad \lambda\in\sigma_\pm^{(1),\mathrm{u}},
\end{equation}
and let us drop some dependencies until the end of this lemma for
notational simplicity.

Let $\lambda\in\inte(\sigma_\mp^{(1),\mathrm{u}})$, then $\phi_\pm,\theta_\pm\in\mathbb{R}$
and $\overline T_\mp=(\overline W\,g_\mp)^{-1}$. Thus, by condition \textbf{(I), (b)} and
\eqref{5.10} for $\lambda\in\inte(\sigma_\mp^{(1)})$ we obtain
\begin{equation}\label{5.29}
\left[\frac{\phi_+\phi_-}{W}\right]=
 \phi_\pm \left[\frac{\phi_\mp}{W}\right] = -g_\mp\phi_\pm
\left(\phi_\mp T_\mp+ \overline\phi_\mp \overline T_\mp\right)=
-g_\mp \theta_\pm \phi_\pm |T_\mp|^2.
\end{equation}
Since $g_\pm\in\mathbb{R}$ for $\lambda\in\inte(\sigma_\mp^{(1),\mathrm{u}})$, \eqref{2.18} implies
\[
\left[\frac{\theta_\mp}{W}\right] = -g_\pm \left[\theta_\mp
T_\pm\right].
\]
The only non-real summand in \eqref{5.11} is the
Cauchy-type integral. The Plemelj formula applied to this
integral gives
\[
\left[\theta_\mp T_\pm\right]= \pm g_\mp\phi_\pm
|T_\mp|^2W(\psi_\pm,\breve\psi_\pm)= g_\mp \phi_\pm|T_\mp|^2
\frac{1}{g_\pm}.
\]
Thus by \eqref{5.29}
\begin{equation}\label{5.301}
\left[\frac{\theta_+\theta_-}{W}\right]=\left[\frac{\phi_+\phi_-}{W}\right]
= - g_\mp \phi_\pm \theta_\pm |T_\mp|^2,\quad
\lambda\in\inte(\sigma_\mp^{(1)}).
\end{equation}
Since $\tilde W\neq 0$ and $s_\mp\neq 0$ for $\lambda\in\inte(\sigma_\mp^{(1)})$, the function
\[
g_\mp \phi_\pm\theta_\pm |T_\mp|^2= - \frac{\delta_\mp^2}{g_\mp}
\frac{\tilde\phi_\pm \tilde\theta_\pm}{|\tilde W|^2}
\]
is bounded on the set under consideration. Finally,
\eqref{5.301} and \eqref{5.28} imply \eqref{5.26}.

(ii) Now we prove, that the function $G(z,x)$ has
removable singularities at the points
$\pa\sigma_-\cup\pa\sigma_+\cup\sigma_d$. Divide this set in four subsets
\begin{equation}\label{5.32}
\Omega_1^\pm=\pa\sigma^{(2)}\cap\inte(\sigma_\mp),\:
\Omega_2^\pm=\pa\sigma_\pm^{(1)}\cap\pa\sigma,\:
\Omega_3=\pa\sigma_-\cap\pa\sigma_+,
\mbox{ and } \Omega_4=\sigma_d.
\end{equation}
In our example we have $\Omega_1^+=\emptyset$, $\Omega_1^-=\{E_5\}$,
$\Omega_2^+=\emptyset$, $\Omega_2^-=\{E_3\}$,
and $\Omega_3=\{E_1,E_2,E_4\}$.

By condition \textbf{II, (b)} all singularities of $G$ are at most isolated poles. Thus it is
sufficient to check that
\begin{equation}\label{5.35}
G(z)=o\left((z -E)^{-1}\right)
\end{equation}
from some direction in the complex plane.

\textbullet\
$\Omega_1$: Consider $E \in \Omega_1^+$ (the case $E \in \Omega_1^-$ is
completely analogous). We will study $\lim_{\lambda\to E}G(\lambda,x)$ as $\la \in \inte(\sigma^{(2)})$
using identity \eqref{5.22}. We have $\phi_-=O(1)$, $g_-=O(1)$, and $\hat W(E)\neq0$.
Moreover, from Lemma~\ref{lem2.1} respectively \textbf{II} we deduce
\[
\phi_+(\lambda)=
\begin{cases}
O(1), & E \notin \hat M_+,\\
O\left(\frac{1}{\sqrt{\lambda-E}}\right), & E \in \hat M_+,
\end{cases}
\qquad
\frac{1}{T_+(\lambda)} =
\begin{cases}
O\left(\frac{1}{\sqrt{\lambda-E}}\right), & E \notin \hat M_+,\\
O(1), & E \in \hat M_+,
\end{cases}
\]
which shows
\[
\theta_-(\lambda)=\frac{\ov{\phi_+(\lambda)}+R_+(\lambda)\phi_+(\lambda)}{T_+(\lambda)}=
O\left(\frac{1}{\sqrt{\lambda-E}}\right),
\]
for $\la \in \sigma^{(2)}$. Inserting this into \eqref{5.22} shows $G(\lambda,x)=O\left(\frac{1}{\sqrt{\lambda-E}}\right)$
and finishes the case $E\in\Omega_1$.

\textbullet\
$\Omega_2$: Suppose now that $E \in \pa \sigma_-^{(1)} \cap\pa \sigma$
(the case $E \in \pa \sigma_+^{(1)} \cap \pa \sigma$ can be treated in the same manner). Now we
cannot use \eqref{5.22}, so we proceed directly from formula \eqref{5.211} estimating the
summands $\frac{\phi_+\phi_-}{W}$ and $\frac{\theta_+\theta_-}{W}$ separately. By
Lemma~\ref{lem2.1} and \textbf{II, (b)} we have
\begin{equation} \label{5.53}
\frac{\phi_+\phi_-}{W}=\frac{\hat \phi_+\hat \phi_-}{\hat W} =
O\left(\frac{1}{\sqrt{z -E}}\right).
\end{equation}
Hence the first term is under control and it remains to investigate the second one.
We investigate the limit of $G$ from the set $\inte(\sigma_-^{(1)})$. By \eqref{5.211}
\begin{equation} \label{5.54}
\frac{\theta_+}{W}= \theta_+T_- g_- = (\ov{\phi_-}+\phi_- R_-) g_-=
\begin{cases}
O(1), & E \in \hat M_-,\\
O\left(\frac{1}{\sqrt{\lambda-E}}\right), & E \notin \hat M_-,
\end{cases}
\;\lambda\in\inte(\sigma_-).
\end{equation}
Next we will estimate $\theta_-$ using \eqref{5.11}. First, let
$E \notin \sigma_v$, that is $\hat W(E)\neq 0$. Using the same
notation, see \eqref{5.123}, as in the proof of Lemma~\ref{lem5.1}
we can split $\theta_-(\lambda)$ according to
\begin{equation} \label{5.55}
\theta_-(\lambda)= \theta_1(\lambda) + \theta_2(\lambda),
\end{equation}
where
\begin{equation} \label{5.57}
\theta_2(\lambda)= \frac{1}{2\pi\I} \int_{\sigma_-^{(1),u}} h_-(\lambda,\xi)
\frac{d\xi}{\xi-\lambda},\qquad \theta_1(\lambda)=\theta_-(\lambda) -
\theta_2(\lambda).
\end{equation}
We see that (cf.\ \eqref{5.142})
\[
\theta_1= O(\zeta_-)= \frac{\hat W}{\hat\delta_-} O(1) =
\begin{cases}
O\left(\frac{1}{\sqrt{\lambda-E}}\right), & E \in \hat M_-,\\
O(1), & E \notin \hat M_-,
\end{cases}
\]
since $E\not\in\sigma_+$. Next, we use \eqref{5.57}, where (cf.\
\eqref{5.13})
\[
h_-(\lambda,\xi)= \frac{\sqrt{\xi-E}}{|\hat W(\xi)|^2}
C(\xi) O(\zeta_-(\lambda))
\]
with $C(\xi)$ some bounded function near $E$.
Hence $\theta_2=O(\zeta_-)$ as well,
which implies together with \eqref{5.54} that
\[
\frac{\theta_+(\lambda)\theta_-(\lambda)}{W(\lambda)}=O\left(\frac{1}{\sqrt{\lambda-E}}\right).
\]
This finishes the case $E\notin\sigma_v$.

Now let $E\in\sigma_v$. As before we have
\begin{equation}\label{5.59}
\theta_1= O(\zeta_-) =
\begin{cases}
O(1), & E \in \hat M_-,\\
O\left(\sqrt{\lambda-E}\right), & E \notin \hat M_-.
\end{cases}
\end{equation}
For the Cauchy-type integral we now have
\[
h_-(\lambda,\xi)= \frac{C(\xi)}{\sqrt{\xi-E}} O(\zeta_-(\lambda))
\]
and \cite[Eq.~(29.8)]{Mu} implies
\begin{align} \label{5.62}
\theta_2(\lambda)=
\begin{cases}
o\left(\frac{1}{\sqrt{\lambda-E}}\right), & E \in \hat M_-,\\
o(1), & E \notin \hat M_-.
\end{cases}
\end{align}
Combining \eqref{5.59}, \eqref{5.62}, and \eqref{5.54} finishes the second case.

\textbullet\
$\Omega_3$: The first step is similar as in the case of $\Omega_2$.
In particular, \eqref{5.53} is valid for $E\in\Omega_3$ and
$z\in\mathbb{C}$ close to $E$. Estimate \eqref{5.54} is valid for
$\lambda\in\inte(\sigma_-)$. The only difference being that $\zeta_-$  in
estimate for $\theta_-$ has an additional square root singularity
since $E\in\pa\sigma_+$. That is, instead of \eqref{5.59} and
\eqref{5.62} we now have
\[
\theta_-(z)= \zeta_-(z) (C + o(1))=
\begin{cases}
O\left(\frac{1}{z -E}\right), & E \in \hat M_-,\\
O\left(\frac{1}{\sqrt{z -E}}\right), & E \notin \hat
M_-,
\end{cases}\quad z\in\mathbb{C}.
\]
However, this is not
good enough unless we can show $C=0$, in which case the big $O$ will
turn into a small $o$ and we are done. It is sufficient to show
$C=0$ from one direction, say $\lambda\in\sigma_+$, which can be done using
the scattering relations for $\theta_-$ as follows.

If $E\in\sigma_v$
this follows directly from
\[
\theta_-(\lambda)=\frac{\overline{\phi_+} +
R_+(\lambda)\phi_+(\lambda)}{T_+(\lambda)}=
O\left(\frac{\hat
W(\lambda)\hat\phi_+(\lambda)}{\hat\delta_-(\lambda)\sqrt{\lambda-E}}\right)=
\begin{cases}
o\left(\frac{1}{\lambda-E}\right), & E \in \hat M_-,\\
o\left(\frac{1}{\sqrt{\lambda-E}}\right), & E \notin \hat
M_-.
\end{cases}
\]
Otherwise, if $E\notin \sigma_v$, then we have two representations
\begin{align}\label{dop5}
\theta_-(\lambda)
&= \frac{1}{T_+(\lambda)} \left(
(\ov{\phi_+(\lambda)}-\phi_+(\lambda)) + \phi_+(\lambda)(R_+(\lambda)+1)
\right),\quad \lambda\notin \hat M_+,\\
\label{dop6}
\theta_-(\lambda)
&= \frac{1}{T_+(\lambda)} \left(
(\ov{\phi_+(\lambda)}+\phi_+(\lambda)) + \phi_+(\lambda)(R_+(\lambda)-1)
\right),\quad \la \in \hat M_+.
\end{align}
For \eqref{dop5} we use
$\ov{\phi_+(\lambda)}-\phi_+(\lambda)=o(1)$ by Lemma~\ref{lem2.1} (ii) and
$R_+(\lambda)+1=o(1)$ by condition \textbf{III}. Now one checks
\[
\left|\frac{1}{T_+(\lambda)}\right|+\left|\frac{\phi_+(z)}{T_+(z)}\right| =
\begin{cases}
O\left(\frac{1}{z -E}\right), & E \in \hat M_-,\\
O\left(\frac{1}{\sqrt{z -E}}\right), & E \notin \hat M_-,
\end{cases}
\quad\lambda\notin\hat M_+,
\]
which shows $G(z,x)=o(\frac{1}{z -E})$ for $\lambda\notin\hat M_+$.
For \eqref{dop6} we use $\ov{\phi_+(\lambda)}+\phi_+(\lambda)=O(1)$ and
$R_+(\lambda)-1=o(1)$ as $\lambda\in\hat M_+$. Since in this case
$\frac{1}{T_+(z)}=O\left(\frac{1}{\sqrt{z -E}}\right)$  and
\[
\frac{\phi_+(z)}{T_+(z)}=
\begin{cases}
O\left(\frac{1}{z -E}\right), & E \in \hat M_-,\\
O\left(\frac{1}{\sqrt{z -E}}\right), & E \notin \hat M_-,
\end{cases}
\quad\lambda\notin\hat M_+,
\]
this implies $G(z,x)=o(\frac{1}{z -E})$ for $\lambda\in\hat M_+$
as required. This finishes the case $E\in \Omega_3$.

\textbullet\
$\Omega_4$: Finally we have to check, that the singularities at the points of the
discrete spectrum are also removable. Since $\tilde W(z)$ has simple zeros at points
$\lambda_k$, it suffices to check that
\begin{equation}\label{5.71}
\tilde\theta_+(\lambda_k,x)\tilde\theta_-(\lambda_k,x) =
\tilde\phi_-(\lambda_k,x)\tilde\phi_+(\lambda_k,x).
\end{equation}
Lemma~\ref{lem5.1} shows, that
$\tilde\theta_\mp:=\delta_\mp\theta_\mp$, defined by \eqref{5.11}, are
continuous at the points  $\breve M_\pm$. Since $(\delta_\mp
T_\pm^{-1})(\lambda_k)=0$ and $(\delta_\mp T_\pm^{-1}\breve\psi_\pm)(\lambda_k)=0$, then
the only the last summand in \eqref{5.11} is non-zero. Computing the limit of
this summand at $\lambda_k$ and using \eqref{2.10}, \eqref{2.17} we obtain
\begin{equation}\label{5.77}
\tilde\theta_\mp(\lambda_k) =\frac{d\tilde
W(\lambda_k)}{d\lambda}(\gamma_\pm)^2\tilde\phi_\pm(\lambda_k),
\end{equation}
which together with \eqref{2.11} implies \eqref{5.71}.
\end{proof}

Lemma \ref{lem5.2} implies, that $G(z,x)$ is an entire function for fixed $x$.
Since in addition $G(z,x)\to 0$ as $z \to\infty$, Liouville's theorem implies
$G(z,x)\equiv 0$ for every $x$. Therefore we have the equalities
\begin{equation}\label{5.72}
\phi_+(z,x)\phi_-(z,x)=\theta_+(z,x)\theta_-(z,x)
\end{equation}
and
\begin{equation}\label{5.73}
\phi_\pm(\lambda,x)\overline{\theta_\pm(\lambda,x)} =
\overline{\phi_\pm(\lambda,x)} \theta_\pm(\lambda,x),\quad \lambda\in\sigma^{(2)}.
\end{equation}
It remains to show that $\theta_\pm(z,x)=\phi_\pm(z,x)$. This is equivalent to showing
that
\[
p(z,x):= \frac{\theta_+(z,x)}{\phi_+(z,x)} =
\frac{\phi_-(z,x)}{\theta_-(z,x)}.
\]
is equal to one. This function is well defined as long as $\hat\phi_+(z,x)\ne0$.
If $\hat\phi_+(z,x)=0$ this is still true as long as $\hat\theta_-(z,x)$ has
no singularity (which is the case for $z\not\in\pa\sigma_-$ by Lemma~\ref{lem5.1})
and does not vanish. But for $z\not\in\pa\sigma_-$ the case $\hat\phi_+(z,x)=\hat\theta_-(z,x)=0$
implies $\hat W(z)=0$, that is, $z\in\sigma_d \cup\pa\sigma_+$. Hence we will
avoid such cases and suppose $x\notin X$, where
\[
X:=\bigcup_{\lambda\in\sigma_d\cup\pa\sigma_-\cup\pa\sigma_+}\{x\mid\hat\phi_+(\lambda,x)=0\}.
\]
Fix $x\in\mathbb{R}\setminus X$. Our aim is to show, that
\begin{equation}\label{5.75}
p(z,x)\equiv 1,\quad z\in\mathbb{C},\ x\notin X.
\end{equation}
By Lemma~\ref{lem2.1} and Corollary~\ref{col7} the functions $\tilde\phi_\pm(z,x)$
and $\tilde\theta_\pm(z,x)$ are holomorphic in $\mathbb{C}\setminus \sigma_\pm$
and hence $p(z,x)$ is holomorphic on $\mathbb{C}\setminus \sigma^{(2)}$ with
continuous limits up to the boundary away from $\pa\sigma^{(2)}$. By \eqref{5.73} the
limits from different sides match up and so $p(z,x)$ is even holomorphic on
$\mathbb{C}\setminus \pa\sigma^{(2)}$. Furthermore, arguing as before, one sees
\[
\frac{\theta_+(z,x)}{\phi_+(z,x)} =
\frac{\hat\theta_+(z,x)}{\hat\phi_+(z,x)} =
O\biggl(\frac{1}{\sqrt{z-E}}\biggr),
\quad E\in\pa\sigma_+,
\]
that is, the singularities near $E\in\pa\sigma^{(2)}$ are removable and
so $p(z,x)$ is entire with respect to $z$ for all $x\notin X$.
Finally, $p(z,x) \to 1$ as $z\to \infty$ by item (iii) of Lemma~\ref{lem1.1} respectively
Lemma~\ref{lem5.1}. In summary, \eqref{5.75} holds, that is,
\begin{equation}\label{5.76}
\theta_\pm(z,x)\equiv\phi_\pm(z,x)
\end{equation}
for all $x\notin X$. But since the set $X$ is discrete, this even holds for all $x\in\mathbb{R}$
by continuity with respect to $x$.

Finally, \eqref{5.10} shows that $\tilde H_\pm\theta_\mp(z,x) = z\theta_\mp(z,x)$, that is,
$\tilde q_+(x)\equiv q_-(x)$. Moreover, \eqref{5.10} and \eqref{5.76} imply, that
\[
T_\mp(\lambda)\phi_\pm(\lambda,x)=\overline{\phi_\mp(\lambda,x)} +
R_\mp(\lambda)\phi_\mp(\lambda,x),
\]
and from \eqref{5.77}, \eqref{5.76}, \eqref{2.11}, and \eqref{D.3} we conclude
that \eqref{2.14} is valid. Therefore the data $\mathcal{S}$ are the scattering data
for the Schr\"odinger operator with the potential $q(x)=\tilde q_-(x)=\tilde q_+(x)$.
This finishes the proof of Theorem~\ref{theor2}.
\end{proof}

Finally, observe that the second moment in condition \textbf{IV} from Lemma~\ref{lem4.1} can be
replaced by any higher moment. In fact, introduce\\[3mm]
\textbf{IV*}.
{\it The function $F_\pm(x,y)$ is continuously differentiable with respect to both variables and
there exists real-valued function $q_\pm(x)$, $x\in\mathbb{R}$, with
\[
\pm\int_a^{\pm\infty}(1+|x|^n) |q_\pm(x)| d x<\infty,  \quad\mbox{for all } a\in\mathbb{R},
\]
such that
\begin{align}\nn
& |F_\pm(x,y)|\leq C_\pm(x)
Q_\pm\left(x+y\right),\\ \nn
& \left|\frac{\pa}{\pa x} F_\pm (x,y)\right|\leq
C_\pm(x)\left(\left|q_\pm\left(\frac{x+y}{2}\right)\right|+Q_\pm(x+y)\right),\\ \nn
& \pm\int_a^{\pm\infty}\,\left|\frac{d}{d x}F_\pm (x,x)\right|(1+ |x|^n)\,d x<\infty,
\end{align}
where $n=2,3,4,\dots$,
\[
Q_\pm(x):=\pm\int_{\frac{x}{2}}^{\pm\infty}\left|q_\pm(t)\right|d t,
\]
and $C_\pm(x)>0$ is a continuous function, which decreases monotonically as $x\to\pm\infty$.}

Then, proceeding literally as in Theorem~\ref{theor2}, we obtain

\begin{theorem}\label{theor3}
Let the scattering data ${\mathcal S}$, defined as in \eqref{4.6},
satisfy conditions \emph{\textbf{I}\textendash\textbf{III}}, and \emph{\textbf{IV*}}. Then each
of the GLM equations \eqref{ME} has a unique solution $K_\pm(x,y)$.
The functions $\tilde q_\pm(x)$ defined by \eqref{5.1}, coincide:
\[
\tilde q_-(x)=\tilde q_+(x)=q(x)
\]
and satisfy
\begin{equation} \label{6.16}
\pm \int_a^{\pm \infty}| q(x) - p_\pm(x)|
(1+|x|^n) d x <\infty ,\quad a\in\mathbb{R}.
\end{equation}
\end{theorem}

\section{The Korteweg--de Vries equation with steplike finite-gap initial data}
\label{KdV}

In this final section we will outline in what way the inverse
scattering transform method can be used to study the initial-value
problem for  the Korteweg--de Vries equation with steplike finite-gap
initial data. Note that the Cauchy problem for steplike constant initial data
is studied in \cite{C,Kap}. In the case of a constant background on
the right half-axis and periodic finite-gap background on the left one, this
problem is partly considered in \cite{Er1}. The existence of the solution of the KdV
equation with general finite-gap steplike potential as initial data
seems to be an open problem and is not  a subject of the present
paper. Here we only discuss a possible approach.

Let $q(x)$ be a smooth function, satisfying condition \eqref{6.16} for $n=5$ together
with its derivatives:
\begin{equation} \label{6.26}
\pm \int_a^{\pm \infty}| q^{(k)}(x) - p_\pm^{(k)}(x)| (1+|x|^5) d x <\infty,
\quad a\in\mathbb{R}, \quad k=0,1,\dots,7.
\end{equation}
Here $p_\pm(x)$ are some finite-gap potentials.

Consider the initial value problem
\begin{align}\label{kdv}
& \frac{\partial u}{\partial t} - 6 u\ \frac{\partial u}{\partial x} + \frac{\partial^3 u}{\partial x^3}=0,\\
\label{id} & u(x,0)=q(x),\quad x\in\mathbb{R}.
\end{align}
One can ask to find a unique smooth solution of this problem in the
domain $|t|<T$, satisfying  conditions
\begin{align}\label{6.1}
& \sup_{|t|<T} \pm\int_0^{\pm\infty} (1+|x|^2) |u(x,t) - u_\pm(x,t)| d x <\infty,\\ \label{6.21}
& \sup_{|t|<T} \pm\int_0^{\pm\infty}
(1+|x|)\left|\frac{\pa^k u(x,t)}{\pa x^k} - \frac{\pa^k
u_\pm(x,t)}{\pa x^k}\right| d x <\infty,\quad k=1,2,3,
\end{align}
where the functions $u_\pm(x,t)$ are the finite-gap solutions of
equation \eqref{kdv} with initial data $u_\pm(x,0)=p_\pm(x)$. We
will proceed by a standard scheme
\[
u(x,0)\ \leadsto\  \mathcal S(0)\ \leadsto\ \mathcal S(t)\ \leadsto\ u(x,t).
\]
The Lax pair, associated with the KdV equation has the form
\begin{align} \label{Aop}
P(t) &= -4\frac{\pa^3}{\pa x^3} + 6u(x,t) \frac{\pa}{\pa x} + 3 u_x(x,t),\\ \label{Hop}
H(t) &=-\frac{\pa^2}{\pa x^2} + u(x,t).
\end{align}
Equation \eqref{kdv} is then equivalent to equation $\pa_t H = [H,P]$.
Let $H_\pm(t)$, $P_\pm(t)$ be Lax pairs, corresponding to our backgrounds.
Following the scheme proposed in \cite{EMT1}, we will obtain the time-dependent
GLM equation:

\begin{lemma}\label{lem12}
Let $\psi_\pm(z,x,t)$ be the Weyl solutions satisfying (see, e.g.,
\cite{GH})
\begin{equation}
H_\pm(t)\psi = z \psi,\qquad P_\pm(t) \psi =
\frac{\pa}{\pa t} \psi,\qquad \psi_\pm(z,0,0)=1.
\end{equation}
Then the
inverse scattering problem (the GLM) equation has the form
\begin{equation}\label{ME1}
K_\pm(x,y,t) + F_\pm(x,y,t) \pm \int_x^{\pm\infty}
K_\pm(x,s,t) F_\pm(s,y,t)d s =0, \quad \pm y>\pm x,
\end{equation}
where
\begin{align}\label{6.2}
F_\pm(x,y,t)
&=
\oint_{\sigma_\pm} R_\pm(\lambda,0)\,\psi_\pm(\lambda,x,t) \psi_\pm(\lambda,y,t)
d\rho_\pm(\lambda,0) \\
&\quad
+\,\int_{\sigma_\mp^{(1),\mathrm{u}}} |T_\mp(\lambda,0)|^2 \psi_\pm(\lambda,x,t)
\psi_\pm(\lambda,y,t)d\rho_\mp(\lambda,0)\notag\\
&\quad
+ \sum_{k=1}^p \gamma_k^\pm(0)^2 \tilde\psi_\pm(\lambda_k,x,t) \tilde\psi_\pm(\lambda_k,y,t).\notag
\end{align}
\end{lemma}

To prove that the problem \eqref{kdv}--\eqref{id} has a solution in
the class \eqref{6.26}--\eqref{6.21} one has to check, that
time-dependent scattering data satisfy conditions of Theorem~\ref{theor2}.
Conditions \textbf{I}, \textbf{II}, \textbf{III (a)}  can be  verified directly from the
time evolution for scattering data
\begin{align*}
R_\pm(\lambda,t)
&=R_\pm(\lambda,0)\exp\bigl(\alpha_\pm(\lambda,t)
-\overline{\alpha_\pm(\lambda,t)}\bigr),\\
T_\pm(\lambda,t)
&=T_\pm(\lambda,0)\exp\bigl(\overline{\alpha_\pm(\lambda,t)} -\alpha_\mp(\lambda,t)\bigr),\\
\gamma_k^\pm(t)^2
&=\gamma_k^\pm(0)^2 \exp (2\alpha_\pm(\lambda_k,t)),
\end{align*}
where (cf.~\cite{Er1,F4})
\[
\alpha_\pm(\lambda,t)=\int_0^t\left(2(u_\pm(0,s) + 2\lambda)m_\pm(\lambda,s) -
\frac{\partial u_\pm(0,s)}{\partial x}\right)ds
\]
and $m_\pm(\lambda,t)$ are time-dependent
Weyl functions, corresponding to background operators $H_\pm(t)$.
Note that condition \textbf{III, (b)} means, that when $E\notin\sigma_v$
and $\mu_j^\pm(t)$ gets close to $E$, then $R_\pm(E,t)$ changes its
sign. This effect was explained in \cite{EMT1}.

To check condition \textbf{IV} one has to take into account the
structure of the Weyl solutions $\psi_\pm(\lambda,x,t)$. As is known
(see e.g., \cite{BBEIM}), they admit a representation
$\psi_\pm(\lambda,x,t)=\exp(\pm \I k_\pm(\lambda)x) f_\pm(\lambda,x,t)$, where
$k_\pm(\lambda)$ are the quasi-momentum maps. To estimate the parts
$F_{d,\pm}(x,y,t)$ and $F_{h,\pm}(x,y,t)$ of the kernel $F_\pm$
(cf.\ \eqref{4.14}, \eqref{4.21}, \eqref{4.24} it is sufficient to
use the Herglotz property of the quasi-momentum ($\Im(k_\pm(\lambda))>0$
as $\lambda\in\mathbb{R}\setminus \sigma_\pm$) and the fact, that the
functions $f_\pm(\lambda,x,t)$ are bounded as $x\in\mathbb{R}$, $|t|<T$
and $\lambda\notin M_\pm(0)\cup \hat M_\pm(0)$, $\la>
\inf\{E_0^+,E_0^-,\lambda_1\} - 1$. To estimate $F_{r,\pm}(x,y,t)$ we
integrate the first summand in \eqref{6.2} twice by parts with
respect to the quasi-momentum variable $k_\pm$ and  prove, that the
boundary terms vanish. This approach fails only in points of the set
$(\pa\sigma_-^{(1)}\cup\pa\sigma_+^{(1)})\cap\pa\sigma^{(2)}$ (the points of
type $E_5$ in our example). Here one has to use the approach,
developed in \cite[Proposition~2.7]{Kap}. This way one arrives at the
estimates
\[
|F_\pm(x,y,t)|\leq\frac{C(t)}{|x+y|^3},\quad \left|\frac{\pa F_\pm(x,y,t)}{\pa x}\right|\leq
\frac{C(t)}{|x+y|^4},\quad x,y\to\pm\infty,
\]
that justify condition \eqref{6.1}.

\subsection*{Acknowledgments}
I.E. gratefully acknowledge the extraordinary hospitality both of the Department of Mathematics
at Paris 7 during a one month stay in the summer of 2005 and of the Faculty of
Mathematics at the University of Vienna during extended stays 2006--2007, where parts
of this paper were written. In addition, we are also grateful to Johanna Michor for discussions on this
topic and to K.\ Grunert for pointing out errors in a previous version of this article.

\appendix
\section{Properties of the transformation operators and estimates for the GLM kernel}
\label{secTOP}

In this appendix we derive and thoroughly investigate the integral equations
for the kernels $K_\pm(x,y)$ of the transformation operators. In addition,
we will obtain the necessary estimates for them and their
derivatives. This will allow us to simplify the necessary and sufficient
conditions on the functions $F_\pm(x,y)$ (in comparison with
\cite{F2}) and to solve the scattering problem in the prescribed
class of perturbations \eqref{0.2}.

Most of the results from this section are in essence known
or follow as in the case of a constant background (see,
e.g., \cite{F1,F2}, and in the discrete case \cite{EBM1,EMT,Khan}).
We included them here (with proofs) to make our
presentation self-contained.

Let $\psi_\pm(z,x)$ be the background Weyl solution \eqref{1.5}. Set
\begin{equation}\label{A.1}
J_\pm(z,x,y)=\frac{\psi_\pm(z,y)\breve\psi_\pm(z,x) -
\psi_\pm(z,x)\breve\psi_\pm(z,y)}{W(\psi_\pm(z),\breve\psi_\pm(z))}
\end{equation}
and
\begin{equation}\label{A.2}
q_\pm(x)=q(x) - p_\pm(x).
\end{equation}
Then the Jost solutions \eqref{2.2} satisfy the integral equation
\begin{equation}\label{A.3}
\phi_\pm(z,x) =\psi_\pm(z,x) -
\int_x^{\pm\infty}J_\pm(z,x,y)q_\pm(y)\phi_\pm(z,y)d y.
\end{equation}
If we substitute formula \eqref{2.2} into this equation, multiply with
$\breve\psi_\pm(z,s)g_\pm(z)$, and integrate over the set $\sipmul$,
using the inverse Fourier transform \eqref{5.4}, and taking into
account that $K_\pm(x,y)=0$, $\pm x >\pm y$, we obtain
\begin{equation}\label{A.4}
K_\pm(x,s) +\int_x^{\pm\infty} dy\,
q_\pm(y)\oint_{\sigma_\pm}J_\pm(\lambda,x,y)\psi_\pm(\lambda,y)\breve\psi_\pm(\lambda,s)
d\rho_\pm(\lambda)\qquad\qquad
\end{equation}
\[
\qquad\pm\int_x^{\pm\infty}d y\, q_\pm(y)\int_y^{\pm\infty} dt\,
K_\pm(y,t)\oint_{\sigma_\pm}J_\pm(\lambda,x,y)\psi_\pm(\lambda,t)\breve\psi_\pm(\lambda,s)d\rho_\pm(\lambda)= 0.
\]
Set
\begin{equation}\label{A.5}
\Gamma_\pm(x,y,t,s)=\mp\oint_{\sigma_\pm}
\psi_\pm(\lambda,x)\breve\psi_\pm(\lambda,y)\psi_\pm(\lambda,t)\breve\psi_\pm(\lambda,s) g_\pm(\lambda)d\rho_\pm(\lambda),
\end{equation}
where the integral has to be understood as a principal value.

Then substituting \eqref{1.12}, \eqref{1.62}, \eqref{A.1}, and \eqref{A.5} into \eqref{A.4} we obtain
\begin{align}\label{A.6}
&K_\pm(x,s) +\int_x^{\pm\infty}\left(\Gamma_\pm(x,y,y,s) -
\Gamma_\pm(y,x,y,s)\right)q_\pm(y)\,d y\\
&\pm\int_x^{\pm\infty}d y \,q_\pm(y)\int_y^{\pm\infty}K_\pm(y,t)\left(\Gamma_\pm(x,y,t,s)
- \Gamma_\pm(y,x,t,s)\right)d t=0.\notag
\end{align}
Consider now the function in \eqref{A.5}. From \eqref{1.10} it follows that
\begin{equation}\label{A.7}
\overline{\Gamma_\pm(x,y,t,s)}=-\Gamma_\pm(y,x,s,t).
\end{equation}
Our plan is to evaluate the integral in \eqref{A.5} using the Jordan lemma.
The only poles of the integrand in \eqref{A.5} are at the band edges and hence
we introduce
\begin{align}\label{A.8}
f_\pm(E,x,y) =& \lim_{z\to E} \Big(\prod_{j=1}^{r_\pm} (z-\mu_j^\pm) \Big)\psi_\pm(z,x)\breve\psi_\pm(z,y)
\end{align}
and
\begin{equation}\label{A.9}
D_\pm(x,y,t,s)=\pm \frac{1}{4}\sum_{E\in\pa\sigma_\pm}
\frac{f_\pm(E,x,y)f_\pm(E,t,s)}{\frac{d}{d z}{P}_\pm(E)}, \quad
P_\pm(z)=\prod_{j=0}^{2r_\pm} (z-E_j^\pm).
\end{equation}
Note that $D_\pm(x,y,t,s)$ is a continuous and bounded function with respect
to all variables.

Now suppose $\pm(x-y+t-s)>0$ and take a closed contour consisting of a large circular arc
together with some parts wrapping around the spectrum $\sigma_\pm$ inside this
arc at a small distance from the spectrum. Due to \eqref{1.62}, \eqref{1.8}, and (iii) of
Lemma~\ref{lem1.1}, it follows that
\[
g_\pm(z)^2 \psi_\pm(z,x)\breve\psi_\pm(z,y)\psi_\pm(z,t)\breve\psi_\pm(z,s)=
O\left(\frac{1}{z}\right) \E^{\pm\I\sqrt{z}(x-y+t-s)}
\]
as $z\to\infty$.
In fact this holds on the entire circle since the neighborhood of the positive real axis can be handled
as above. Hence one can apply Jordan's lemma to conclude that the contribution of the circle
vanishes as its radius tends to infinity. Shrinking the loops the integral converges to
\begin{equation}\label{A.10}
\Gamma_\pm(x,y,t,s)=D_\pm(x,y,t,s),\quad\mbox{for}\quad
\pm(x-y+t-s)>0.
\end{equation}
Note that $f_\pm(E,x,y)$ are real, and
$f_\pm(E,x,y)=f_\pm(E,y,x)$. Thus, the $D_\pm(x,y,t,s)$ is also real,
\begin{equation}\label{A.11}
D_\pm(x,y,t,s)=D_\pm(y,x,t,s).
\end{equation}
Now let $\pm(x-y+t-s)<0$, that is, $\pm(y-x+s-t)>0$.
Then \eqref{A.7}, \eqref{A.10}, and \eqref{A.11} imply
\[
\Gamma_\pm(x,y,t,s)=-\overline{D_\pm(x,y,t,s)}=-D_\pm(x,y,t,s),\quad \pm(x-y+t-s)<0.
\]
Therefore,
\begin{equation}\label{A.12}
\Gamma_\pm(x,y,t,s)=D_\pm(x,y,t,s)\sign(\pm(x-y+t-s)).
\end{equation}
Property \eqref{A.11} implies that the domain, where
in the first integrand in \eqref{A.6} does not vanish, is
\begin{equation}\label{A.13}
\sign(\pm(x-s))=-\sign(\pm(2y-x-s)),\quad \pm s>\pm x.
\end{equation}
In the second integral the domain of integration is
\begin{align}\label{A.121}
&\sign(\pm(x-y+t-s))=-\sign(\pm(y-x+t-s)),\\
&\text{ with}\quad\pm s>\pm x,\quad\pm t>\pm y>\pm x.\notag
\end{align}
Solving \eqref{A.13} and \eqref{A.121} we arrive at the following result.

\begin{lemma}\label{lemA.1}
The kernels $K_\pm(x,s)$ of the transformation operators satisfy
the integral equation
\begin{align}\nn
K_\pm(x,s)
&= - 2\int_{\frac{x+s}{2}}^{\pm\infty}q_\pm(y)D_\pm(x,y,y,s)d y\\\label{A.15}
&\quad
 \mp 2\int_x^{\pm \infty}d y\int_{s\pm x\mp y}^{s\pm y\mp
x}D_\pm(x,y,t,s)K_\pm(y,t)q_\pm(y)\,d t,\quad\pm s >\pm x,
\end{align}
where $D_\pm$ are defined by \eqref{A.9}.
\end{lemma}

Set $s=x$ in \eqref{A.15}. Then the second summand vanishes,
because we have our integration inside the domain $\pm t<\pm y$, where
$K_\pm(y,t)=0$. Thus
\begin{equation}\label{A.16}
K_\pm(x,x)=- 2\int_x^{\pm\infty}q_\pm(y)D_{\pm}(x,y,y,x)\,d y.
\end{equation}
But, as is well-known (see, e.g, \cite[Eq.~(1.84)]{GH} or \cite[Chapter~8]{L})
\begin{equation}\label{A.17}
\psi_\pm(z,y)\breve\psi_\pm(z,y)\prod_{j=1}^{r_\pm} (z-\mu_j^\pm)=
\prod_{j=1}^{r_\pm} (z-\mu_j^\pm(y)),
\end{equation}
where $\mu_j^\pm(y)$ are the Dirichlet eigenvalues corresponding to the
base point $x=y$ (rather than $x=0$). Combining
\eqref{A.5} and \eqref{A.10} we obtain
\begin{align}
D_\pm(x,y,y,x)
&= D_\pm(x,x,y,y)\notag\\ \label{A.18}
& = \pm\frac{1}{4}\sum_{E\in\pa\sigma_\pm}\Res_E
\frac{\prod_{j=1}^{r_\pm}\left(z -\mu_j^\pm(x)\right)\left( z
-\mu_j^\pm(y)\right)}{(z -E_0^\pm)\prod_{j=1}^{r_\pm} \left((z -
E_{2j-1}^\pm)(z -E_{2j}^\pm)\right)}.
\end{align}
The integrand in \eqref{A.18} is meromorphic in $\mathbb{C}$, thus, by
the Cauchy theorem, we can compute the residue at infinity and obtain
\[
D_\pm(x,y,y,x)=-\lim_{z\to\infty}\pm\frac{1}{4}\frac{z}{z -E_0^\pm}
\prod_{j=1}^{r_\pm}\frac{  \left(z -\mu_j^\pm(y)\right)\left( z
-\mu_j^\pm(x)\right)}{ (z -E_{2j-1}^\pm)(z -
E_{2j}^\pm)}=\mp\frac{1}{4}.
\]
From \eqref{A.16} we conclude that
\begin{equation}\label{A.19}
K_\pm(x,x)=\pm\frac{1}{2}
\int_x^{\pm\infty}(q(t) -p_\pm(t))d t.
\end{equation}
This formula justifies formula \eqref{5.1} under the condition, that  the transformation
operators kernels are differentiable.

\begin{lemma}\label{lemA.2}
Let
\begin{equation}\label{A.20}
Q_\pm(x):=\pm\int_{\frac{x}{2}}^{\pm\infty}\left|q_\pm(t)\right|\,d t,\quad
q_\pm(x)=q(x) - p_\pm(x).
\end{equation}
Then $K_\pm(x,y)$ has first order partial derivatives with respect to both variables.
Moreover, for $\pm y\geq\pm x$ the following estimates are valid
\begin{align}\label{A.21}
\left|K_\pm(x,y)\right|
&\leq C_\pm(x) Q_\pm(x+y),\\
\label{A.22}
\left|\frac{\partial K_\pm(x,y)}{\partial x}\right| +
\left|\frac{\partial K_\pm(x,y)}{\partial y}\right|
&\leq
C_\pm(x)\,\left(\left|q_\pm\left(\frac{x+y}{2}\right)\right|
+ Q_\pm\left(x+y\right)\right),
\end{align}
where $C_\pm(x)$ are positive continuous
functions for $x\in\mathbb{R}$ which decrease as $x\to \pm\infty$
and depend on the corresponding background data and on the first
moment of the perturbation.
\end{lemma}

\begin{proof}
We restrict our considerations to the ``$+$" case only and omit ``$+$" in what follows.
We will follow the scheme of the proof of \cite[Lemmas~3.1.1, 3.1.2]{M}.
Introduce the following change of variables in \eqref{A.15}:
\begin{equation}\label{A.28}
y+t=:2\alpha,\ \ t-y=:2\beta,\ \
x+s=:2u,\ \ s-x=:2v,
\end{equation}
then from \eqref{A.15} we obtain (see \cite[Lemma~3.1.1]{M})
\begin{align}\label{A.23}
H(u,v)
&= -2\int_u^\infty q(s) D_1(u,v,s)\,d s\\
&\quad
- 4\int_u^\infty
d\alpha\int_0^v q(\alpha-\beta)
D_2(u,v,\alpha,\beta)H(\alpha,\beta)\,d\beta,\notag
\end{align}
where we put
\begin{equation}\label{A.24}
\begin{split}
&H(u,v)=K_+(u-v,u + v),\quad
D_1(u,v,s)=D_+(u-v,s,s,u+v),\\
&D_2(u,v,\alpha,\beta)=D_+(u-v,\alpha-\beta,\alpha+\beta,u+v).
\end{split}
\end{equation}
Functions $D_1$ and $D_2$ are bounded uniformly with respect to
all their variables. Put $C=2\max\{\max_{u,v,s}|D_1|,
\max_{u,v,\alpha, \beta}|D_2|\}$ and apply the method of
successive approximations (see \cite[Lemma~3.1.1]{M}). We arrive
at the estimate
\begin{equation}\label{A.25}
|H(u,v)|\leq \tilde C(u-v)\int_u^\infty|\tilde q(x)|\,d x,
\end{equation}
with
\begin{equation}\label{A.40}
\tilde C(u)=C\exp\left(C\int_{2u}^\infty Q(2t)d t\right),\quad C>0,
\end{equation}
from which \eqref{A.21} follows. To obtain \eqref{A.22}, observe that the
first partial derivatives of $D_1$ and $D_2$ exist (see
\eqref{A.8}, \eqref{A.9} and \eqref{A.24}) and are bounded with respect to
all variables. Thus,
\begin{align}\label{A.26}
&\frac{\partial H(u,v)}{\partial
u}-2 q(u)D_1(u,v,u)=\notag\\
&= - 4\int_0^v q(u - \beta)
D_2(u,v,u,\beta)H(u,\beta)\,d\beta
-2\int_u^\infty q(s) \frac{\pa D_1(u,v,s)}{\pa u}\,d s\notag\\
&\quad
-4\int_u^\infty d\alpha\int_0^v q(\alpha - \beta) \frac{\pa
D_2(u,v,\alpha,\beta)}{\pa u}H(\alpha,\beta)\,d\beta,
\end{align}
\begin{align}\label{A.27}
&\frac{\partial H(u,v)}{\partial v}=\notag\\
&=
-2\left(\int_u^\infty q(s) \frac{\pa D_1(u,v,s)}{\pa
v}\,d s - 2\int_u^\infty q(\alpha - v)
D_2(u,v,\alpha,v)H(\alpha,v)\,d\alpha\right.\notag\\
&\qquad\quad
-\left. 2 \int_u^\infty d\alpha\int_0^v q(\alpha - \beta)
\frac{\pa D_2(u,v,\alpha,\beta)}{\pa
v}H(\alpha,\beta)\,d\beta\right).
\end{align}
The function $Q(u)=\int_u^\infty| q(x)|\,dx$ is positive,
monotonically decreasing, and satisfies
$Q(\,\cdot\,)\in L^1(a,\infty)$, $a\in\mathbb{R}$. Since
$|D_i|$, $|\frac{\pa D_i}{\pa u}|$, $|\frac{\pa D_i}{\pa v}|\leq
C_1$, $i=1,2$, \eqref{A.25} applied to
\eqref{A.26} and \eqref{A.27} implies
\[
\left|\frac{\partial H(u,v)}{\partial u}-2 q(u)D_1(u,v,u)\right| +
\left|\frac{\partial H(u,v)}{\partial v}\right|\leq \tilde C(u-v) Q(2u),
\]
where the function $\tilde C(u)$ is of the same type as \eqref{A.40}, with a different
positive constant $C_2$ depending on the background data. From this,
\eqref{A.28}, and \eqref{A.24}, the estimate \eqref{A.22} follows.
\end{proof}

With the help of this lemma we can now derive several estimates for the GLM equation.

\begin{lemma}\label{lemA.3}
The kernel $F_\pm(x,y)$
of the GLM equation \eqref{ME} has first order derivatives with
respect to each variable. Furthermore, for $\pm y>\pm x$ it satisfies
\begin{equation}\label{A.29}
\left|F_\pm(x,y)\right|\leq \hat C_\pm(x)\,Q_\pm(x+y),
\end{equation}
\begin{equation}\label{A.30}
\left|\frac{\partial F_\pm(x,y)}{\partial x}\right| + \left|\frac{\partial
F_\pm(x,y)}{\partial y}\right|\leq \hat
C_\pm(x)\,\left(\left|q_\pm\left(\frac{x+y}{2}\right)\right| +
Q_\pm\left(x+y\right)\right),
\end{equation}
where the functions $q_\pm(x)$
and $Q_\pm(x)$ are defined in \eqref{A.20}. Here $\hat C_\pm(x)$ are positive continuous
functions which decrease as $x\to\pm\infty$. Moreover,
\begin{equation}\label{A.32}
\pm\int_a^{\pm\infty}(1+x^2)\left|\frac{d F_\pm(x,x)}{d x}\right|<\infty,\quad \forall
a\in\mathbb{R}.
\end{equation}
\end{lemma}

\begin{proof}
Again we restrict our considerations to the ``$+$" case only and omit ``$+$" in
what follows. Set $Q_1(u)=\int_u^\infty Q(t)d t$. Due to
condition \eqref{0.2}, the  functions $Q(x)$ and $Q_1(x)$ satisfy
\begin{equation}\label{A.31}
\int_a^\infty Q_1(t)d t<\infty, \quad \int_a^\infty Q(t)(1+|t|)d t<\infty.
\end{equation}
Observe also, that the kernel $F(x,y)$ of the GLM equation \eqref{ME}
is symmetric: $F(x,y)=F(y,x)$. From \eqref{ME} and \eqref{A.21} we see, that
\[
|F(x,y)|\leq \tilde C(x)\left(Q(x+y) +\int_x^\infty Q(x+t)|F(t,y)|d t\right).
\]
Since $Q_1(x+t)>Q_1(2x)$, and $\tilde C(t)<\tilde C(x)$ as $x<t$, then
Gronwall's inequality implies \eqref{A.29} with
\[
\hat C(x)
= C_1\tilde C(x)\exp\bigl(C_1\tilde C(x) Q_1(2x)\bigr),\quad
C_1>0.
\]
Differentiating \eqref{ME} with respect to $x$ and $y$ implies
\begin{align}\label{A.33}
&
\left|F_x(x,y)\right|\leq\left|K_x(x,y)\right| +
|K(x,x)F(x,y)| +\int_x^\infty \left|
K_x(x,t)F(t,y)\right|d t,\\
\label{A.34}
&
F_y(x,y) + K_y(x,y) +\int_x^\infty K(x,t)
F_y(t,y)d t=0.
\end{align}
The functions $Q(x)$, $Q_1(x)$, $\hat C(x)$, $\tilde C(x)$ are monotonously decreasing
and positive. Furthermore,
\[
\int_x^\infty\left(
\left|q_\pm\left(\frac{x+t}{2}\right)\right|
+Q(x+t)\right)Q(t+y)d t\leq(Q(2x)
+Q_1(2x))Q(x+y),
\]
and hence the estimate \eqref{A.30} for $F_x$ follows (with some other positive continuous
decreasing function $\hat C(x)$) from \eqref{A.29}, \eqref{A.21},
\eqref{A.22} and \eqref{A.33}. The same estimate for $F_y$ can be obtained from \eqref{A.34}
and Lemma~\ref{lemA.2} by using the method of successive approximations.

It remains to prove \eqref{A.32}. To this end consider \eqref{ME}
for $y=x$ and differentiate it with respect to $x$:
\[
\frac{d F(x,x)}{d x} + \frac{d K(x,x)}{d x} - K(x,x) F(x,x)
+\int_x^\infty \left(K_x(x,t) F(t,x) + K(x,t)
F_y(t,x)\right)d t=0.
\]
Formula \eqref{A.19} implies \eqref{5.101}.  Next, by \eqref{A.29} and \eqref{A.21}, we have
\[
|K(x,x)F(x,x)|\leq\tilde C(a)\hat C(a) Q^2(2x)\ \mbox{for} \ x>a,
\]
where $\int_a^\infty(1+x^2) Q^2(2x)d x<\infty$. Moreover,
by \eqref{A.30} and \eqref{A.22},
\[
\left|K_x^\prime (x,t) F(t,x)\right| + \left|K(x,t) F_y^\prime
(t,x)\right|\leq 4\tilde C(a)
\hat C(a)\Bigl\lbrace\Bigl\lvert q\Bigl(\frac{x+t}{2}\Bigr)\Bigr\rvert Q(x+t) +
Q^2(x+t)\Bigr\rbrace,
\]
and together with the estimates
\begin{align*}
& \int_a^\infty d x\,x^2\int_x^\infty Q^2(x+t)d t\leq\int_a^\infty
|x| Q(2x)d x\ \sup_{x\geq a}\int_x^\infty|x+t| Q(x+t)d t<\infty,\\
& \int_a^\infty x^2\int_x^\infty\Bigl\lvert q\Bigl(\frac{x+t}{2}\Bigr)\Bigr\rvert
Q(x+t)d t\leq\\
& \qquad
\leq\int_a^\infty Q(2x)d x \ \sup_{x\geq a} \int_x^\infty
\Bigl\lvert q\Bigl(\frac{x+t}{2}\Bigr)\Bigr\rvert(1 +(x+t)^2)d t<\infty
\end{align*}
we arrive at \eqref{A.32}.
\end{proof}

\begin{remark}\label{rem6}
Note that the results of this lemma are in some sense invertible.
Namely, if we start with properties \eqref{4.3}--\eqref{4.31} of $F_\pm(x,y)$, then,
using \eqref{ME} and the same considerations as in Lemma~\ref{lemA.3},
we obtain \eqref{5.100}, \eqref{A.22}, and \eqref{5.101}.
\end{remark}

\section{Proof of Lemma \ref{lem2.2}}
\label{appB}

We introduce the local parameter $\tau=\sqrt{z -E}$ in a small
vicinity of each point $E\in\pa\sigma_\pm$ and set $\dot{y}(z,x)
= \frac{\pa}{\pa \tau} y(z,x)$. Since $\frac{d z}{d\tau}(E)=0$,
for every solution $y(z,x)$ of the equation \eqref{1.1}, its derivative
$\dot{y}(E,x)$ is also a solution of \eqref{1.1}. In particular,
the Wronskian $W(y(E),\dot{y}(E))$ is independent of $x$.

For each $x\in\mathbb{R}$ in a small neighborhood of a point $E\in\pa\sigma_\pm$ introduce the function
\begin{equation}\label{B.0}
\hat\psi_{\pm,E}(z,x) =
\begin{cases}
\psi_\pm(z,x),&E\,\in\pa\sigma_\pm\setminus \hat M_\pm,\\
\tau\, \psi_\pm(z,x),&E\,\in \hat M_\pm.
\end{cases}
\end{equation}

\begin{lemma}\label{lem1.2}
Let the function $\hat\psi_{\pm,E}(z,x)$ be defined by formula \eqref{B.0}.
Then
\begin{equation}\label{1.15}
W\Bigl(\hat\psi_{\pm,E}(E), \frac{\pa}{\pa \tau}\hat{\psi}_{\pm,E}(E)\Bigr) =
\pm \lim_{z\to E} \frac{\alpha\,\tau^\alpha}{2g_\pm(z)},
\end{equation}
where $\alpha=-1$ if $ E\,\in\pa\sigma_\pm\setminus \hat M_\pm$ and
$\alpha=1$ if $E\,\in \hat M_\pm$.
\end{lemma}

\begin{proof}
We begin by recalling (see, e.g., \cite[Eq.~(1.73)]{GH}) that the Weyl $m$-functions can be written as
\begin{equation}
m_\pm(z) = \frac{G_\pm(z) \pm \sqrt{P_\pm(z)}}{F_\pm(z)},
\end{equation}
where $G_\pm(z)=\frac{1}{2}F_\pm'(z,0)$ and $F_\pm(z)=F_\pm(z,0)$ with
\begin{equation}
P_\pm(z)=\prod_{j=0}^{2r_\pm} (z-E_j^\pm), \qquad
F_\pm(z,x)=\prod_{j=1}^{r_\pm} (z-\mu_j^\pm(x)).
\end{equation}
Furthermore, observe that $\mp \sqrt{P_\pm(z)}= \tau f_\pm(z)$, where $f_\pm(z)$ is
holomorphic near $z=E$.

Then in the case $E\in\pa\sigma_\pm\setminus\hat M_\pm$ (where $\hat\psi_{\pm,E}(\lambda,x)= \psi_\pm(\lambda,x)$)
we have $F_\pm(E)\ne 0$. Then, from \eqref{1.5}, we have $\psi_\pm(E,0)=1$,
$\psi_\pm^\prime(E,0)=m_\pm(E)$, $\dot \psi_\pm(E,0)=0$,
\[
\dot\psi_\pm^\prime(E,0)=\dot m_\pm(E) = \frac{f_\pm(E)}{F_\pm(E)}
= \mp \lim_{z\to E} \frac{\sqrt{P_\pm(z)}}{\tau F_\pm(z)},
\]
and the first claim follows.

In the second case we have $E\in\hat M_\pm$ (where $\hat\psi_{\pm,E}(\lambda,x)= \tau \psi_\pm(\lambda,x)$)
we have $F_\pm(E)= 0$ and hence $F_\pm(z)=\tau^2 \tilde{F}_\pm(z)$,
$G_\pm(z)=\tau^2 \tilde{G}_\pm(z)$, where $\tilde{F}_\pm(z)$ and $\tilde{G}_\pm(z)$ are holomorphic
near $z=E$ with $\tilde{F}_\pm(E)\ne 0$ and $\tilde{G}_\pm(E) \ne 0$. Hence we
have
\begin{alignat*}{2}
&\hat\psi_{\pm,E}(E,0)=0,&\qquad&
\hat\psi_{\pm,E}^\prime(E,0) =  \frac{f_\pm(E)}{\tilde{F}_\pm(E)},\\
&
\frac{\pa}{\pa\tau}\hat\psi_{\pm, E}(E,0) =1,&&
\frac{\pa}{\pa \tau} \hat\psi_{\pm,E}^\prime(E,0) =  \frac{\tilde{G}_\pm(E)}{\tilde{F}_\pm(E)}
\end{alignat*}
and the second claim follows.
\end{proof}

From \eqref{1.88}, \eqref{1.4} and \eqref{1.5} we observe the following property of the Floquet--Weyl solutions

\begin{remark}\label{rem3}
Let $E\in\hat M_\pm$ and $x\in\mathbb{R}$.
Then $\overline{\hat \psi_{\pm,E}(E,x)}=-\hat\psi_{\pm,E}(E,x)$ if $E$ is a left band edge from $\sigma_\pm$
and $\overline{\hat \psi_{\pm,E}(E,x)}= \hat\psi_{\pm,E}(E,x)$ if $E$ is a right band edge.
\end{remark}

It is straightforward to check that this property is also inherited by the Jost solutions.
Again we abbreviate
\begin{equation}\label{2.361}
\hat\phi_{\pm,E}(\lambda,x) =
\begin{cases}
\phi_\pm(\lambda,x),&E\,\in\pa\sigma_\pm\setminus \hat M_\pm,\\
\tau\,\phi_\pm(\lambda,x),&E\,\in \hat
M_\pm.
\end{cases}
\end{equation}

\begin{lemma}\label{lem2.1a}
We have $\overline{\phi_\pm(E,x)}=\phi_\pm(E,x)$ for $E\in\pa\sigma_\pm\setminus \hat M_\pm$.
If $E\in\hat M_\pm$, then by \eqref{2.361}
$\overline{\hat \phi_{\pm,E}(E,x)}=-\hat\phi_{\pm,E}(E,x)$ when
$E$ is a left edge of a band from $\sigma_\pm$ and $\overline{\hat
\phi_{\pm,E}(E,x)}=\hat\phi_{\pm,E}(E,x)$ when $E$ is a right
edge of a band.
\end{lemma}

\begin{lemma}\label{lem2.2}
The function $\hat W(z)$ is continuous on the set $\mathbb{C}\setminus\sigma$
up to the boundary $\siu\cup\sil$. It can have zeros on the set
$\pa\sigma\cup(\pa\sigma_+^{(1)}\cap\pa\sigma_-^{(1)})$ and does
not vanish at the other points of the set $\sigma$.  If $\hat W(E)=0$
as $E\in\pa\sigma\cup(\pa\sigma_+^{(1)}\cap\pa\sigma_-^{(1)})$,
then $\hat W(z) =  \sqrt{z -E} (C(E)+o(1))$, $C(E)\ne 0$.
\end{lemma}

\begin{proof}
Continuity of $\hat W$ up to the boundary follows from the corresponding
property of $\hat \phi_\pm(z,x)$. We begin with the investigation of
the possible zeros.

Let $\lambda_0\in\inte(\sigma^{(2)}):= \sigma^{(2)}\setminus\pa\sigma^{(2)}$ and
suppose $W(\lambda_0) = 0$. Then $\phi_+(\lambda_0,x)=c\,\phi_-(\lambda_0,x)$,
$\overline{\phi_+(\lambda_0,x)}=\bar c\, \overline{\phi_-(\lambda_0,x)}$,
i.e. $W(\phi_+,\overline {\phi_+})=|c|^2W(\phi_-,\overline
{\phi_-})$. But by \eqref{1.62} and \eqref{1.12} we have $\sign
g_+(\lambda_0)= -\sign g_-(\lambda_0)$, contradicting \eqref{1.881}.

Let $\lambda_0\in\inte(\sigma^{(1)}_\pm)$ and $\tilde W(\lambda_0) =0$. The
point $\lambda_0$ can coincide with a pole $\mu\in M_\mp$. But
$\phi_\pm(\lambda_0,x)$ and $\overline{\phi_\pm(\lambda_0,x)}$ are linearly
independent and bounded, and $\tilde\phi_\mp(x,\lambda_0)\in\mathbb{R}$
(see \eqref{2.12}). If $W(\lambda_0)=0$, then
$\tilde\phi_\mp=c_1^{\pm}\,\phi_\pm=c_2^\pm\,\overline{\phi_\pm}$,
that implies $W(\phi_\pm,\overline{\phi_\pm})(\lambda_0)=0$,
which is impossible.

In the general mutual location of the background spectra one can
meet the case when
$\lambda_0=E\in\pa\sigma^{(2)}\cap\inte(\sigma_\pm)\subset\inte(\sigma)$ (that is
a point like $E_5$ in our example). If $\hat W(E)=0$, then
$W( \phi_\pm, \hat\phi_{\mp,E})(E)=0$, where
$\hat\phi_{\mp,E}$ is defined by \eqref{2.361}. But according to
Lemma~\ref{lem2.1a}, the values of
$\hat\phi_{\mp,E}(E,\,\cdot\,)$ are either pure real or pure imaginary,
therefore $W(\overline{\phi_\pm},\hat\phi_{\mp,E})(E)=0$, that is, $\overline{\phi_\pm(E,x)}$ and $\phi_\pm(E,x)$ are linearly dependent, which is impossible at
inner points of the set $\sigma_\pm^{(1)}$.

In summary, $\hat W(\lambda)\neq 0$ for $\lambda\in\inte(\sigma)\setminus(\pa\sigma_-^{(1)}\cap\pa\sigma_+^{(1)})$ which finishes the first part and it remains to investigate the
order of zeros.

Let $E\in\pa\sigma\cup(\pa\sigma_+^{(1)}\cap\pa\sigma_-^{(1)})$ (these are
the points of type $E_1$, $E_3$ and $E_4$ from our example). The
function $\hat W(\lambda)$ is continuously differentiable with respect
to the local parameter $\tau$. Since at $E$ we have
$\frac{d}{d\tau}(\delta_+ \delta_-)(E)=0$, the function
$W(\hat\phi_+,\hat\phi_-)$ has the same order of zero at $E$ as
$\hat W(\lambda)$. But if $\hat\delta(E)\neq 0$, then
$\frac{d}{d\tau}\hat \delta_\pm(E)=0$ and if
$\hat\delta_-(E)=\hat\delta_+(E)= 0$, then
$\frac{d}{d\tau}(\tau^{-2}\,\hat \delta_+\,\hat\delta_-)(E)=0$.
Therefore $\frac{d}{d\tau}\hat W(E) = 0$ if and only if
$\frac{d}{d\tau}W(\hat\phi_{+,E},\hat\phi_{-,E}) = 0$.

To simplify notations, we will just write $\hat\phi_\pm:=\hat\phi_{\pm,E}$
until the end of this proof. Again we have consider all possible cases for the
mutual location of the Dirichlet eigenvalues.

First let $E\in\pa\sigma^{(2)}\cap\pa\sigma$ (a point of type $E_2$ in our example)
and let $E\notin(\hat M_+\cup\hat M_-)$. Then $\hat W(E)=0$ if and only if
(see \eqref{2.8}) $W(E)=W(\phi_+,\phi_-)=0$, that is, $\phi_\pm(E,\,\cdot\,)
=c_\pm\,\phi_\mp(E,\,\cdot\,)$, $c_-c_+=1$, $c_-,c_+\in\mathbb{R}$.
The derivative of the Jost solution with respect to $\tau$ is again
a solution of equation \eqref{2.1}. Therefore, by Lemma~\ref{lem1.2},
\begin{align}\label{D5}
\dot W(E)
&=W(\dot \phi_+,\phi_-) - W(\dot
\phi_-,\phi_+)=c_-\, W(\dot \phi_+,\phi_+) -
c_+\,W(\dot \phi_-,\phi_-)\\
&= c_-\,W(\dot \psi_+,\psi_+) - c_+\,W(\dot \psi_-,\psi_-)
=-(c_+\,d_- +c_-\,d_+),\notag
\end{align}
where
\[
d_\pm=
\lim_{\lambda\to E} \frac{\I}{2 g_\pm(\lambda)\sqrt{\lambda-E}}.
\]
We see from \eqref{1.881} that $d_\pm\in \I\mathbb{R}_+\setminus\{0\}$ if $E$
is a left edge of $\sigma$ and $d_\pm\in \mathbb{R}_+\setminus\{0\}$ if $E$ is a right edge of
$\sigma$. Since $\sign c_-=\sign c_+$, this finishes the case
$E\in\pa\sigma^{(2)}\cap\pa\sigma\setminus(\hat M_+\cup\hat M_-)$.

The same arguments are valid for
$E\in\pa\sigma^{(2)}\cap\pa\sigma\cap\hat M_+\cap\hat M_-$ and $\hat W(E)=0$. Then $\hat
\phi_\pm(E,\,\cdot\,)=c_\pm\,\hat\phi_\mp(E,\,\cdot\,)$ and using Lemma~\ref{lem2.1a}
we conclude that $c_-,\,c_+\in\mathbb{R}$, $\sign c_-=\sign c_+$. By Lemma~\ref{lem1.2}
\[
\frac{d}{d\tau}W(\hat\phi_+,\hat\phi_-)(E) =
c_+ \lim_{\lambda\to E}\frac{\sqrt{\lambda-E}}{2g_-(\lambda)}
 +c_- \lim_{\lambda\to E}\frac{\sqrt{\lambda-E}}{2g_+(\lambda)}\neq 0.
\]

The cases, when one of the Jost solutions is bounded in the edge of spectrum,
and the another one is unbounded, are more subtle. For example, let $\hat W(E)=0$
for $E\in(\hat M_+\cap\pa\sigma^{(2)}\cap\pa\sigma)\setminus(\hat M_-\cap\hat M_+)$
and let $E$ be a right band edge (like the point $E_2$
of or example when $\mu_+=E_2$ and $\mu_-\neq E_2$). Then by
Lemma~\ref{lem2.1a} $\hat\phi_+(x,E)\in\mathbb{R}$ and
$\hat\phi_-(x,E)=\phi_-(x,E)\in\mathbb{R}$. Also
$\frac{d}{d\tau}\hat W(E)\neq 0$ if and only if
$\frac{d}{d\tau}W(\hat\phi_+,\,\phi_-)(E)\neq 0$.
Therefore we have $\hat\phi_+=c_+\phi_-$, $\phi_-=c_-\hat\phi_+$
and $\sign c_+=\sign c_-$. Thus, by Lemma~\ref{lem1.2}
$\frac{d}{d\tau}W(\hat\phi_+,\,\phi_-)(E)= c_- d_+ - c_+\,d_-$, where
\[
d_+ = \lim_{\lambda\to E}\frac{\I\sqrt{\lambda-E}} {2\I g_+(\lambda)}\in\mathbb{R}_-\setminus \{0\}, \quad
d_- = \lim_{\lambda\to E} \frac{-1}{\I\sqrt{\lambda-E}\, 2\I g_-(\lambda)}\in\mathbb{R}_+\setminus\{0\}.
\]
The case, when $E\in(\hat M_-\cap\pa\sigma^{(2)}\cap\pa\sigma)
\setminus(\hat M_+\cap\hat M_-)$ and $E$ is a right
band edge is analogous.

Now, let $E\in(\hat M_+\cap\pa\sigma^{(2)}\cap\pa\sigma)\setminus\hat M_-$
and let $E$ be a left band edge. Then by Lemma~\ref{lem2.1a}
$\hat\phi_+\in \I\mathbb{R}$ and $c_+, c_-\in \I\mathbb{R}$,
$c_+\,c_-=1$, that is, $\sign ( \I c_+)=-\sign (\I c_-)$.
Furthermore, $\sign\sqrt{\lambda-E}>0$ since $\la>E$, and $\sign
g_+(\lambda)=\sign g_-(\lambda)$. Therefore,
\[
\frac{d}{d\tau} W(\hat\phi_+,\phi_-)(E)=\I c_-\lim_{\lambda\to E}\frac{
\sqrt{\lambda-E}} {g_+(\lambda)}- \I\ c_+\lim_{\lambda\to E}\frac{
1}{\sqrt{\lambda-E}\,g_-(\lambda)} \neq 0.
\]

Unlike the case of one and the same background $p_+(x)=p_-(x)$,
where $\hat W(\lambda)\neq 0$ for $\lambda\in\inte(\sigma)$, we could have
$\hat W(\lambda)=0$ for $\lambda\in\inte(\sigma)$ in our steplike situation.
The points under consideration are points of the
set $\pa\sigma^{(1)}_-\cap\pa\sigma^{(1)}_+$. This case can be treated in
the same way as the case $E\in\pa\sigma^{(2)}\cap\pa\sigma$.
Namely, if $E\notin (\hat M_-\cup\hat M_+)$ then observe that one
of the two summands in \eqref{D5} is  real and the other one is
imaginary. The same is valid for the case $E\in (\hat M_-\cap\hat
M_+)$. Now let $E\in\pa\sigma^{(1)}_-\cap\pa\sigma^{(1)}_+\cap\hat
M_+\setminus(\hat M_+\cap\hat M_-)$ and let $E$ be a
right band edge of $\sigma_-^{(1)}$ (i.e.\ a left band edge of
$\sigma_+^{(1)}$). Then by Lemma \ref{lem2.1a} $\hat\phi_+(x,E)\in
\I\mathbb{R}$ and $\hat\phi_-(x,E)=\phi_-(x,E)\in \mathbb{R}$.
Moreover, $\hat W=0$ if and only if $W(\hat\phi_+,\hat\phi_-)=0$ with
the same order of zero. Therefore,  $c_+,c_-\in \I\mathbb{R}$ and by
\eqref{1.15} and \eqref{1.12}
\begin{align*}
\frac{d}{d\tau}W(\hat\phi_+,\hat\phi_-)
&=
c_- W(\hat{\dot{\psi}}_+,\hat\psi_+) - c_+ W( \hat{\dot{\psi}}_-,\hat\psi_-)\\
&= \frac{1}{2\pi}\Bigl(\I c_-\lim_{\lambda\to E}\,\frac{\tau}{g_+} +\I
c_+\lim_{\lambda\to E} \frac{1}{\tau g_-} \Bigr)\neq 0
\end{align*}
because the first summand is imaginary and the second one is real.
All other combinations can be treated similarly.

Finally, consider the case $E\in\pa\sigma^{(1)}_\pm\cap\pa\sigma.$ Let, for
example, $E\in\pa\sigma_+^{(1)}\cap\pa\sigma$ (the point $E_4$ in our example).
Since $E\in\mathbb{R}\setminus\sigma_-$ in this case, we have $\hat\delta_-(E)\neq 0$
and one has to study zero of the Wronskian $W(\hat\phi_+,\tilde \phi_-)$. But by \eqref{1.6}
and \eqref{2.2} $W\bigl(\tilde\phi_-,\frac{d}{d\tau}\tilde\phi_-\bigr)= W\bigl(\tilde\psi_-,\frac{d}{d\tau}\tilde\psi_-\bigr)=0$.
Therefore,
\[
\frac{d}{d\tau}W(\hat\phi_+,\tilde\phi_-)=c_- W\Bigl(\frac{d}{d\tau} \hat\phi_+,\hat\phi_+\Bigr)\neq 0
\]
by Lemma~\ref{lem1.2}.
\end{proof}


\end{document}